\documentclass[10pt,twoside,english]{amsart}
\advance\oddsidemargin by -1.0cm
\advance\evensidemargin by -1.0cm
\advance\topmargin by -1.0cm 
\usepackage{amssymb}
\usepackage{babel}
\usepackage{amstext}
\usepackage{diagrams}
\usepackage{amscd}   
\usepackage{epsfig}  
\usepackage{rotating}
\def\haken{\mathbin{\hbox to 6pt{%
                    \vrule height0.4pt width5pt depth0pt                  
                    \kern-.4pt                  
                    \vrule height6pt width0.4pt depth0pt\hss}}}     
      
\let\mathg\mathfrak
\theoremstyle{plain}
\newtheorem{cor}{Corollary}[section]
\newtheorem{lem}{Lemma}[section]
\newtheorem{thm}{Theorem}[section]
\newtheorem{prop}{Proposition}[section]
\theoremstyle{definition}

\newtheorem{NB}{Remark}[section]
\newtheorem{dfn}{Definition}[section]
\newtheorem{assume}{Assumption}[section]
\newtheorem*{thank}{Thanks}

\newcommand{\bdm}{\begin{displaymath}}
\newcommand{\edm}{\end{displaymath}}
\newcommand{\ba}[1]{\begin{array}{#1}}
\newcommand{\ea}{\end{array}}
\newcommand{\bea}[1][]{\begin{eqnarray#1}}
\newcommand{\eea}[1][]{\end{eqnarray#1}}
\newcommand{\btab}{\begin{tabular}}
\newcommand{\etab}{\end{tabular}}

\newcommand{\x}{\times}

\newcommand{\ox}{\otimes}

\newcommand{\ra}{\rightarrow}

\newcommand{\lmapsto}{\longmapsto}

\newcommand{\lan}{\left\langle}
\newcommand{\ran}{\right\rangle}

\newcommand{\hut}{\wedge}

\newcommand{\tr}{\ensuremath{\mathrm{tr}}}

\newcommand{\del}{\partial}

\newcommand{\ad}{\ensuremath{\mathrm{ad}}}

\newcommand{\R}{\ensuremath{\mathbb{R}}}

\newcommand{\M}{\ensuremath{\mathcal{M}}}

\newcommand{\vphi}{\ensuremath{\varphi}}    
\newcommand{\vrho}{\ensuremath{\varrho}}

\newcommand{\Jacm}{\ensuremath{\mathrm{Jac}_{\mathg{m}}}}
\newcommand{\Jach}{\ensuremath{\mathrm{Jac}_{\mathg{h}}}}
\newcommand{\Jac}{\ensuremath{\mathrm{Jac}}}
\newcommand{\Ric}{\ensuremath{\mathrm{Ric}}}
\newcommand{\Scal}{\ensuremath{\mathrm{Scal}}}

\newcommand{\Ad}{\ensuremath{\mathrm{Ad}\,}}
\newcommand{\LC}{\ensuremath{\mathrm{LC}}}
\newcommand{\Adtilde}{\ensuremath{\widetilde{\mathrm{A}\,}\!\mathrm{d}\,}}
\newcommand{\adtilde}{\ensuremath{\widetilde{\mathrm{a}\,}\!\mathrm{d}\,}}
\newcommand{\Ctilde}{\ensuremath{\widetilde{C\,}\!}}
\renewcommand{\ad}{\ensuremath{\mathrm{ad}\,}}

\newcommand{\diag}{\ensuremath{\mathrm{diag}}}

\newcommand{\Cl}{\ensuremath{\mathcal{C}}}

\newcommand{\GL}{\ensuremath{\mathrm{GL}}}

\newcommand{\SU}{\ensuremath{\mathrm{SU}}}

\newcommand{\so}{\ensuremath{\mathg{so}}}
\newcommand{\SO}{\ensuremath{\mathrm{SO}}}
\newcommand{\Spin}{\ensuremath{\mathrm{Spin}}}
\newcommand{\spin}{\ensuremath{\mathg{spin}}}

\newcommand{\g}{\ensuremath{\mathfrak{g}}}
\newcommand{\h}{\ensuremath{\mathfrak{h}}}
\newcommand{\m}{\ensuremath{\mathfrak{m}}}

\begin{document}
\setcounter{equation}{0}
%
%
\thispagestyle{empty}
%
\date{\today}
\title[Connections on naturally reductive spaces]{Connections  on naturally 
reductive spaces, their Dirac operator and homogeneous models in string theory}
%
%
%
\author{Ilka Agricola}
\address{\hspace{-5mm} 
{\normalfont\ttfamily agricola@mathematik.hu-berlin.de}\newline
Institut f\"ur Reine Mathematik \newline
Humboldt-Universit\"at zu Berlin\newline
Sitz: WBC Adlershof\newline
D-10099 Berlin\\
Germany}
\thanks{This work was supported by the SFB 288 "Differential geometry
and quantum physics" of the Deutsche Forschungsgemeinschaft.}
\subjclass[2000]{Primary 53 C 27; Secondary 53 C 30}
\begin{abstract}
Given a  reductive homogeneous space $M=G/H$ endowed with a 
naturally reductive metric, we study the one-parameter family of connections
$\nabla^t$ joining the canonical and the Levi-Civita connection ($t=0, 1/2$). 
We show that the Dirac operator $D^t$ corresponding to $t=1/3$ is the so-called
``cubic'' Dirac operator recently introduced  by B.\ Kostant, and
derive the formula for its square for any $t$, thus generalizing the
classical Parthasarathy formula on symmetric spaces. Applications include
the existence of a new $G$-invariant first order differential operator
$\mathcal{D}$ on spinors and an eigenvalue estimate for the first
eigenvalue of $D^{1/3}$. This geometric situation can be used
for constructing Riemannian manifolds which are Ricci flat and admit a
 parallel spinor with respect to some metric connection $\nabla$ whose
torsion $T\neq 0$ is a $3$-form, the geometric model for the common sector
of string theories. We present some results about solutions to the
string equations and give a detailed discussion of some 
$5$-dimensional example.
\end{abstract}
\maketitle
\tableofcontents
\pagestyle{headings}
%
%
%
\section{Introduction}\noindent
This paper proposes a differential geometric approach to some recent
results from  B.\ Kostant on an algebraic object called
"cubic Dirac operator" (\cite{Kostant99}). The key observation is
that one can introduce a metric connection on certain homogeneous spaces whose
torsion  (viewed as a ($0,3$)-tensor) is  $3$-form such that the
associated Dirac operator has Kostant's algebraic object as its symbol.
At the same time, there has been recently a growing interest in
connections with  totally skew symmetric torsion for constructing
models in string theory and supergravity. We show that the mentioned class of 
homogeneous spaces yields interesting candidates for such solutions
and use Dirac operator techniques to prove some vanishing theorems. 

In a first part of this paper, we consider a  reductive 
homogeneous space $M=G/H$ endowed with a 
Riemannian metric  that induces a naturally reductive metric $\lan\ , \ \ran$ 
on $\m$, where we set $\g=\h\oplus\m$. The one-parameter family of 
$G$-invariant connections defined by
 \bdm
 \nabla^t_X Y\ =\ \nabla^0_X Y + t\,[X,Y]_{\m}
 \edm
joins the canonical ($t=0$) and the Levi-Civita ($t=1/2$) connection.
Its torsion $T(X,Y,Z)= (2t-1)\cdot\lan[X,Y]_{\m},Z\ran$ is a $3$-form.
For an orthonormal basis $Z_1,\ldots,Z_n$ of $\m$, it induces the third degree 
element 
 \bdm
 H\ :=\  \frac{3}{2}\sum_{i<j<k}\lan[Z_i,Z_j]_{\m},Z_k \ran Z_i\cdot Z_j
\cdot Z_k
 \edm
inside the Clifford algebra $\Cl(\m)$ of $\m$. The fact that the Dirac operator
associated with the connection $\nabla^t$ may then be written as
 \bdm
 D^t\psi\ =\ \sum_{i} Z_i\cdot Z_i(\psi) + t\cdot H\cdot\psi
 \edm
suggested  the name "cubic Dirac operator" to B.\ Kostant. We will show that
the main achievement in \cite{Kostant99} was to realize that, for
the parameter value $t=1/3$, the square of $D^t$ may be expressed in
a very simple way in terms of Casimir operators and scalars only
(\cite[Thm 2.13]{Kostant99}, \cite[10.18]{Sternberg99}). It is a remarkable
generalization of the well-known Parthasarathy formula for $D^2$ on
symmetric spaces (Theorem~\ref{Kos-Parth-D2-symm} in this article, see
\cite{Parthasarathy72}).
In fact, S.\ Slebarski has already noticed independently that the
parameter value $t=1/3$ has distinguished properties (see Theorem 1 and
the introduction in \cite{SlebarskiI87}). He uses it to prove a "vanishing 
theorem" for the kernel of the twisted Dirac operator, which can be easily
recovered from Kostant's formula (see \cite[Thm 4]{Landweber00}).
Although his articles \cite{SlebarskiI87} and \cite{SlebarskiII87} contain
several attempts to generalize Parthasarathy's formula for $D^2$, none of
them seems to come close to Kostant's results. We shall compute the
general expression for $(D^t)^2$ in Theorem~\ref{K-P-1} and show how it 
can be simplified for this particular parameter value in Theorem~\ref{K-P-2}. 
We emphasize one difference between our work and \cite{Kostant99}. While 
Kostant studies the algebraic action
of $D^{1/3}$ as an element of $\mathcal{U}(\g)\ox\Cl(\m)$ on
$L^2$-functions $G\ra\Delta_{\m}$ (the spinor representation), we restrict
our attention to spinors, i.\,e., $L^2$-sections of the spinor
bundle $S=G\x_{\Adtilde}\Delta_{\m}$. In particular, this implies that
one of the terms in the formula for $(D^t)^2$ (the "diagonally" embedded 
Casimir operator of $\h$) vanishes independently of $t$. An immediate
consequence of Theorem \ref{K-P-1} is the existence
of a new $G$-invariant first order differential operator 
\bdm
\mathcal{D}\psi \ :=\ \sum_{i,j,k} \lan [Z_i,Z_j]_{\m},Z_k \ran 
Z_i\cdot Z_j\cdot Z_k(\psi)
\edm
on spinors (Remark~\ref{D-new-op}) that has no analogue on symmetric spaces.
Furthermore, under some additional hypotheses (the lifted Casimir operator 
$\Omega_{\g}$ has to be non negative) Theorem~\ref{K-P-2} yields  
an eigenvalue estimate, which is discussed in Corollary
\ref{eigenvalue-estimate}.

In the second part of this paper, we use the preceding approach
for studying the string equations on naturally reductive spaces.
Stated in a differential geometric way, one wants to construct
a Riemannian manifold $(M,g)$ with a metric connection $\nabla$
such that its torsion $T\neq 0 $ is a $3$-form and such that there exists at 
least one spinor field $\psi$ satisfying the coupled system
 \bdm
 \Ric^{\nabla}\ =\ 0,\quad \delta(T)\ =\ 0,\quad \nabla\Psi\ =\ 0,\quad
 T\cdot \Psi \ =\ 0\,.
 \edm
The number of preserved supersymmetries depends essentially on the number
of $\nabla$-parallel spinors. For a general background on these equations, 
we refer to the article by A.\ Strominger \cite{Strominger86}, where they 
appeared for the first time. Thus, if one looks for homogeneous solutions, 
the family of connections $\nabla^t$ yields
canonical candidates for the desired connection $\nabla$, and the results
on the associated Dirac operator can be used to discuss the solution
space to these equations. We discuss the significance
of constant spinors  (which do not always exist) in 
Theorem~\ref{constant-fields} and show that the last two string
equations cannot have any solutions at all if the lifted Casimir 
operator $\Omega_{\g}$ is non negative (Theorem~\ref{eq-3-4-no-sol}). 
In order to  discuss the first equation, we present a representation 
theoretical expression for the Ricci tensor of the connection $\nabla^t$,
which generalizes previous results by Wang and Ziller 
(Theorem~\ref{gen-Ric-WZ}).
The article ends with a thourough discussion of an example, namely, the
naturally reductive metrics on the $5$-dimensional Stiefel manifold.

Although we rarely refer to it, this paper is in spirit very close
(and in some sense complementary) to a recent article by Friedrich and 
Ivanov (\cite{Friedrich&I01}). There, the authors study  metric connections 
with totally skew symmetric torsion preserving a given geometry. 
%
\begin{thank}
I am grateful to Thomas Friedrich (Humboldt-Universit\"at zu Berlin) for many
valuable discussions on the topic of this paper. My thanks are also
due to the Erwin-Schr\"odinger Institute in Vienna and the Max-Planck
Institute for Mathematics in the Natural Sciences in Leipzig for their 
hospitality.
\end{thank}

\section{A family of connections on naturally reductive spaces}
\label{fam-conn}\noindent
%
Consider a Riemannian homogeneous space $M=G/H$. We suppose that $M$ is
\emph{reductive}, i.\,e.,  the Lie algebra $\g$ of $G$ may be decomposed
into a vector space direct sum of the Lie algebra $\h$ of $H$ 
and an $\Ad(H)$-invariant subspace $\m$ such that $\g=\h\oplus\m$
and $\Ad(H)\m\subset \m$. We identify $\m$ with $T_0M$ by the map
$X\mapsto X_0^*$, where $X^*$ is the Killing vector field on $M$ generated
by the one parameter group $\exp(tX)$ acting on $M$. We pull back
the Riemannian metric $\lan\ ,\ \ran_0$ on $T_0M$ to an inner product
$\lan\ ,\ \ran$ on $\m$.
Let $\mathrm{Ad}:\, H\ra \SO(\m)$ be the isotropy representation of $M$. By a 
theorem of Wang (\cite[Ch. X, Thm 2.1]{Kobayashi&N2}), there is a one-to-one 
correspondence between the set of $G$-invariant metric affine connections 
and the set of linear mappings $\Lambda_{\m}:\ \m\ra\so(\m)$ such that
 \bdm
 \Lambda_{\m}(hXh^{-1})\ =\ \Ad(h)\Lambda_{\m}(X)\Ad(h)^{-1} \
 \text{ for }X\in\m \text{ and } h\in H\,.
 \edm
Its torsion and curvature are then given for $X,Y\in\m$ by 
(\cite[Ch. X, Prop. 2.3]{Kobayashi&N2})
 \bea[*]
 T(X,Y)& =& \Lambda_{\m}(X)Y-\Lambda_{\m}(Y)X - [X,Y]_{\m},\\
 R(X,Y)& =& [\Lambda_{\m}(X),\Lambda_{\m}(Y)] - \Lambda_{\m}([X,Y]_{\m})
 -\Ad([X,Y]_{\h})\, ,
 \eea[*]
where the Lie bracket is split into its $\m$ and $\h$ part,
$[X,Y]=[X,Y]_{\m}+[X,Y]_{\h}$. 
\begin{lem}\label{antisymmetry}
The $(0,3)$-tensor corresponding to the torsion $(X,Y,Z\in\m)$
 \bdm
 T(X,Y,Z)\ :=\ \lan T(X,Y),Z\ran 
 \edm
is totally skew symmetric if and only if the map $\Lambda_{\m}$
satisfies for all $X,Y,Z\in\m$ the invariance condition
 \bdm
 \lan \Lambda_{\m}(X)Y,Z\ran + \lan \Lambda_{\m}(Z)Y,X\ran\ =\ 
 \lan [X,Y]_{\m},Z\ran + \lan [Z,Y]_{\m},X\ran\,.
 \edm
\end{lem}
\begin{proof}
The antisymmetry of $T(X,Y,Z)$ in $X$ and $Z$ is equivalent to
 \bdm
 \lan \Lambda_{\m}(X)Y,Z\ran + \lan \Lambda_{\m}(Z)Y,X\ran -
 \lan \Lambda_{\m}(Y)X,Z\ran - \lan \Lambda_{\m}(Y)Z,X\ran
 -\lan [X,Y]_{\m},Z\ran - \lan [Z,Y]_{\m},X\ran\ =\ 0\,.
 \edm
The third and fourth term cancel out each other by the assumption
that $\Lambda_{\m}(Y)$ lies in $\so(\m)$, since this means that the
endomorphism $\Lambda_{\m}(Y)$ is skew symmetric with respect to the 
inner product of $\m$.
\end{proof}
\noindent
For a general map $\Lambda_{\m}$, this is  all one
can say. We are interested in the one parameter family of connections
defined by
 \bdm
 \Lambda_{\m}^t (X)Y\ :=\ t\cdot [X,Y]_{\m}\,.
 \edm
It is well known that $t=0$ corresponds to the \emph{canonical connection} 
$\nabla^0$, which, by the Ambrose-Singer theorem, is the unique metric 
connection on $M$ such that its torsion and curvature are parallel,
$\nabla^0 T^0=\nabla^0 R^0=0$. By Lemma~\ref{antisymmetry}, the
torsion of  $\nabla^0$ is a $3$-form if and only if $M$ is
\emph{naturally reductive}.
\begin{dfn}
A homogeneous Riemannian metric on $M$ is said to be \emph{naturally
reductive} (with respect to $G$) if the map $[X,-]_{\m}:\m\ra\m$ is 
skew symmetric,
 \bdm
 \lan [X,Y]_{\m},Z\ran + \lan Y, [X,Z]_{\m}\ran \ =\ 0 \text{ for all }
 X,Y,Z\in\m\,.
 \edm
Note that if $G_1\subset G_2$ are two transitive groups of isometries
of $M$, then the properties of being naturally reductive with respect to
$G_1$ and $G_2$ are independent of each other.
\end{dfn}
\begin{NB}\label{skew-symm-of-Lambda}
Under the assumption that $M$ is naturally reductive, the right-hand side
in the criterion of Lemma~\ref{antisymmetry} vanishes, and the remaining
condition may be restated -- using the skew symmetry of $\Lambda_{\m}(X)$
and $\Lambda_{\m}(Z)$ --
as $\lan Y,\Lambda_{\m}(X)Z+ \Lambda_{\m}(Z)X\ran =0$. Since this equation
has to hold for all $X,Y$ and $Z$ in $\m$, we obtain that the torsion
is a $3$-form if and only if $ \Lambda_{\m}(X)X = 0$  for all $X\in\m$.
\end{NB}
\noindent
If $M$ is naturally reductive, then the torsion of the family $\nabla^t$ of 
connections is given by the simple expression
 \bdm
 T^t(X,Y)\ =\ (2t-1)\, [X,Y]_{\m}\,.
 \edm
One sees that the Levi-Civita connection is attained for $t=1/2$. The 
general formula for the connection $\nabla^t$ is
 \begin{equation}\label{connection}
 \nabla^t_X Y\ =\ \nabla^0_X Y + t\,[X,Y]_{\m}\, .
 \end{equation}
Notice that for a symmetric space, $[\m,\m]\subset\h$, so all connections
of this one-parameter family coincide and are equal to the Levi-Civita 
connection.
\begin{assume}
We will assume that $M=G/H$ is naturally reductive with respect to $G$.
\end{assume}
\noindent
We begin by computing a few characteristic entities for this family of 
connections, which will be needed in the subsequent sections. We start by
recalling a theorem of B.\ Kostant.
\begin{thm}[\cite{Kostant56}]
Suppose $G$ acts effectively on $M=G/H$. If the inner product
$\lan\ , \ \ran$ is naturally reductive with respect to $G$,
then $\tilde{\g}:=\m+[\m,\m]$ is an ideal in $\g$ whose
corresponding subgroup $\tilde{G}\subset G$ is transitive on $M$,
and there exists a unique $\Ad(\tilde{G})$ invariant, symmetric,
non degenerate, bilinear form $Q$ on $\tilde{\g}$ (not necessarily
positive definite) such that 
 \bdm
 Q(\h\cap\tilde{\g},\m)\ =\ 0\ \text{ and }\ Q|_{\m}\ =\ \lan\ , \ \ran\, ,
 \edm
where $\h\cap\tilde{\g}$ will be the isotropy algebra in $\tilde{\g}$.
Conversely, if $G$ is connected, then, for any $\Ad(G)$ invariant,
symmetric, non degenerate, bilinear form $Q$ on $\g$, which is
non degenerate on $\h$ and positive definite on $\m:=\h^{\perp}$, the
metric on $M$ defined by $Q|_{\m}$ is naturally reductive. In this
case, $\g=\tilde{\g}$. \qed
\end{thm}
\begin{assume}\label{assume-transitivity}
We shall assume from now on that $G$ acts transitively on $M$
(thus, $\g=\tilde{\g}$) and use the $\Ad(G)$ invariant extension $Q$ of
the inner product $\lan\ , \ \ran$ as well as its restriction
$Q|_{\h}=:Q_{\h}$ to $\h$ where needed without further comment.
\end{assume}
%
%
\begin{lem}\label{curvature}
%
The curvature of the connection $\nabla^t$ is given by
 \bdm
 R^t(X,Y)Z\ =\ t^2\,[X,[Y,Z]_{\m}]_{\m}+ t^2\,[Y,[Z,X]_{\m}]_{\m}+
 t\,[Z,[X,Y]_{\m}]_{\m} + [Z,[X,Y]_{\h}]\,.
 \edm
If $Z_i,\,\ldots,Z_n$ is an orthonormal basis of $\m$, the Ricci tensor
and the scalar curvature are
 \bea[*]
 \Ric^t(X,Y)& =& \sum_i (t-t^2)\lan [X,Z_i]_{\m},[Y,Z_i]_{\m} \ran +
 Q_{\h}([X,Z_i],[Y,Z_i])\\
 \Scal^t &=& \sum_{i,j} (t-t^2)\lan [Z_i,Z_j]_{\m},[Z_i,Z_j]_{\m} \ran +
 Q_{\h}([Z_i,Z_j],[Z_i,Z_j])\,.
 \eea[*]
\end{lem}
\begin{proof}
The formula for the curvature follows immediately from the general formula 
given before. In particular, it implies
 \bdm
 \lan R^t(X,Z)Z,Y\ran \ =\ (t-t^2)\,\lan [X,Z]_{\m},[Y,Z]_{\m}\ran
 + \lan [Z,[X,Z]_{\h}],Y\ran\,.
 \edm
Using the $\Ad(G)$ invariant extension $Q$ of the inner product
$\lan\ , \ \ran$ and the fact that $\m$ is then perpendicular to $\h$, we 
may rewrite the latter term as
 \bdm
 \lan [Z,[X,Z]_{\h}],Y\ran\ =\  Q( [Z,[X,Z]_{\h}],Y)\ =\ 
 Q([X,Z]_{\h},[Y,Z])\ =\ Q_{\h}([X,Z],[Y,Z])\,.
 \edm
Thus, we obtain
 \bdm
 \lan R^t(X,Z)Z,Y\ran \ =\ (t-t^2)\, Q_{\m}( [X,Z],[Y,Z])
 +Q_{\h}([X,Z],[Y,Z])
 \edm
and the formula for the Ricci tensor by 
$\Ric^t(X,Y)=\sum \lan R^t(X,Z_i)Z_i,Y\ran$. The expression for the
scalar curvature is  obtained by contraction relative to $X$ and $Y$.
\end{proof}
\noindent
At a later stage, we will give a further expression for the Ricci tensor
due to Wang and Ziller (\cite{Wang&Z85}). For the time being, we observe
that the connection with $t=1$ has also special properties, for example,
it has the same Ricci tensor than the canonical connection. This is
why we propose to call it the \emph{anticanonical} connection.
We compute the covariant derivative of the torsion tensor.
\begin{lem}\label{nablaT-2-1}
As a map $\m\x\m\ra\m$, the covariant derivative of $T$ is
 \bdm
 (\nabla^t_Z T^t)(X,Y)\ =\ t(2t-1)\big([X,[Y,Z]_{\m}]_{\m}+
[Y,[Z,X]_{\m}]_{\m}+[Z,[X,Y]_{\m}]_{\m}\big)\,.
 \edm
\end{lem}
\begin{proof}
By definition, the covariant derivative is given by
 \bdm
 (\nabla^t_Z T^t)(X,Y)\ =\ \nabla^t_Z(T^t(X,Y)) - T^t(\nabla^t_Z X,Y)
 - T^t(X,\nabla^t_Z Y)\,.
 \edm
We insert the expression for $\nabla^t$ from equation~(\ref{connection})
 \bea[*]
 (\nabla^t_Z T^t)(X,Y)& =& \nabla^0_Z(T^t(X,Y)) + t[Z,T^t(X,Y)]_{\m}
 -T^t(\nabla^0_Z X + t[Z,X]_{\m},Y)\\
 & &- T^t(X,\nabla^0_Z Y + t[Z,Y]_{\m})\\
 &=&  \nabla^0_Z(T^t(X,Y)) -T^t(\nabla^0_Z X,Y)- T^t(X,\nabla^0_Z Y)\\
 && + t(2t-1)\big([X,[Y,Z]_{\m}]_{\m}+[Y,[Z,X]_{\m}]_{\m}+
 [Z,[X,Y]_{\m}]_{\m}\big)\,.
 \eea[*]
But the third line may be rewritten as $-(2t-1)(\nabla^0_Z T^0)(X,Y)$,
which vanishes by the Ambrose-Singer theorem.
\end{proof}
\noindent
For the first time we encounter here an expression that will play
an important role at different places. Let us define
 \bea[*]
 \Jacm(X,Y,Z)& :=& [X,[Y,Z]_{\m}]_{\m}+[Y,[Z,X]_{\m}]_{\m}+
 [Z,[X,Y]_{\m}]_{\m}\, ,\\
 \Jach(X,Y,Z)& :=& [X,[Y,Z]_{\h}]\,+\,[Y,[Z,X]_{\h}]\,+\, [Z,[X,Y]_{\h}]\,.
 \eea[*]
Notice that the summands of $\Jach(X,Y,Z)$ automatically lie in $\m$
by the assumption that $M$ is reductive. The Jacobi identity for $\g$
implies $\lan \Jacm(X,Y,Z)+\Jach(X,Y,Z),\m \ran =0$. As the
connection $\nabla^t$ is metric, the covariant derivatives of $T$
viewed as a $(0,3)$- resp.\ $(1,2)$-tensor are related by
 \begin{equation}\label{nablaT-3-0}
 (\nabla^t_Z T^t)(X,Y,V)\ =\ \lan (\nabla^t_Z T^t)(X,Y),V\ran
 \ =\ t(2t-1)\lan\Jacm(X,Y,Z),V \ran\,.
 \end{equation}
For completeness, we recall the formula for the exterior derivative of
a differential form in terms of a connection with torsion.
\begin{lem}\label{ext-diff}
If $\omega$ is an $r$-form, then
 \bea[*]
(d\omega)(X_0,\,\ldots,X_r)& =& \sum_{i=0}^r (-1)^i (\nabla_{X_i}\omega)
(X_0,\ldots,\hat{X}_i,\ldots,X_r) \\
&- & \sum_{0\leq i<j\leq r}(-1)^{i+j} 
\omega(T(X_i,X_j),X_0,\ldots,\hat{X}_i,\ldots,\hat{X}_j,\ldots,X_r)\,.
 \eea[*]
\end{lem}
\begin{proof}
%
We start with the general formula for the derivative of an $r$-form $\omega$
(see, for example, \cite[Prop. 3.11]{Kobayashi&N1}), 
\bea[*]
(d\omega)(X_0,\,\ldots,X_r)& =& \sum_{i=0}^r (-1)^i 
X_i(\omega(X_0,\ldots,\hat{X}_i,\ldots,X_r))\\
&+& \sum_{0\leq i<j\leq r}(-1)^{i+j} 
\omega([X_i,X_j],X_0,\ldots,\hat{X}_i,\ldots,\hat{X}_j,\ldots,X_r)\,.
\eea[*]
In the first line, we express every summand in terms of the 
covariant derivative of $\omega$, i.\,e.,
\bea[*]
X_i (\omega(X_0,\ldots,\hat{X}_i,\ldots,,X_r))& =& 
(\nabla_{X_i}\omega)(X_1,\,\ldots,X_r)
+ \omega(\nabla_{X_i}X_0,\,\ldots,\hat{X}_i,\,\ldots,X_r)+\ldots +\\
&+& \omega(X_0,X_1,\,\ldots,\nabla_{X_i}X_r)\,.
\eea[*]
A simple rearrangement of terms together with the definition 
$T(X,Y)=\nabla_X Y-\nabla_Y X - [X,Y]$ of the torsion yields the
result.
\end{proof}
\begin{lem}\label{dT}
The codifferential of the $3$-form $T^t$ vanishes, $\delta T^t=0$,
while its outer derivative is given by 
$dT^t(X,Y,Z,V)=2(2t-1)\cdot\lan \Jacm(X,Y,Z),V \ran$.
\end{lem}
\begin{proof}
For the first claim, one deduces from equation~(\ref{nablaT-3-0}) that
$X\haken \nabla^t_X T^t=0$. Then it follows for the orthonormal basis
$Z_i,\,\ldots,Z_n$ of $\m$ that
 \bdm
 \delta^t T^t\ =\ \sum_{i=1}^n Z_i\haken \nabla^t_{Z_i} T^t\ =\ 0\,.
 \edm
In particular, the divergence of $T$ with respect to $\nabla^t$
coincides with its Riemannian divergence (a more general fact, see 
\cite{Friedrich&I01}), $\delta^t T^t=\delta^{1/2} T^t=0$. Hence
we shall  drop the superscript, as we did in the statement of
the lemma. The second claim follows from Lemma~\ref{ext-diff} by a simple 
algebraic computation.
\end{proof}
\begin{NB}\label{torsion-lie-alg}
We finish this section with a remark about the connection between the
torsion and the Lie algebra structure. If some torsion $3$-form $T$
is given as a fundamental datum and is to be the torsion of the canonical
connection of some space with naturally reductive metric, then the $\m$-part 
of the commutators $[\m,\m]$  may be reconstructed by
 \bdm
 [X,Y]_{\m} \ = \ -\sum_{i} T(X,Y,Z_i)Z_i\,.
 \edm
This formula is fundamental for the point of view taken in the article
\cite{Kostant99} (formula 1.23). The full Lie algebra structure of $\g$ can 
now be viewed as consisting of the torsion $3$-form, the isotropy 
representation and the subalgebra structure of $\h$, with some compatibility 
condition resulting from the Jacobi identity. This point of view will be useful
in the last section, where we will study examples.
\end{NB}
\section{The Dirac operator of the family of connections $\nabla^t$}
\label{Dirac}
\subsection{General remarks and formal self adjointness}
\noindent
Assume that there exists a homogeneous spin structure on $M$, i.\ e.,
a lift $\Adtilde:\ H\ra\Spin(\m)$ of the isotropy representation such 
that the diagram
%
\begin{diagram}
  &                       &  \Spin(\m)\\
  &\ruTo^{\Adtilde} & \uTo_{\lambda}\\
H & \rTo^{\Ad} & \SO(\m)\\
\end{diagram}
%
commutes, where $\lambda$ denotes the spin covering. Moreover, we  
denote by $\adtilde$ the corresponding lift into $\spin(\m)$ of
the differential $\ad:\h\ra\so(\m)$ of $\Ad$. Let
$\kappa:\, \Spin(\m)\ra \GL(\Delta_{\m})$ be the spin representation,
and identify sections of the spinor bundle $S=G\x_{\Adtilde}\Delta_{\m}$
with functions $\psi:\ G\ra\Delta_{\m}$ satisfying
 \bdm 
 \psi(gh)\ =\ \kappa(\Adtilde(h^{-1}))\psi(g)\,.
 \edm
For any $G$ invariant
connection defined by a map $\Lambda_{\m}: \m\ra\so(\m)$, we
consider its lift $\tilde{\Lambda}_{\m}: \m\ra\spin(\m)$, which is
given by $\tilde{\Lambda}_{\m}:= d\lambda^{-1}\circ\Lambda_{\m}$.
Then the the covariant derivative on spinors may be expressed as
(\cite[Lemma 2]{Ikeda75}) 
 \begin{equation}\label{gen-cov-der}
 \nabla_Z\psi \ =\  Z(\psi)+\tilde{\Lambda}_{\m}(Z) \psi
 \end{equation}
and thus the  Dirac operator associated with this connection 
has the form
 \begin{equation}\label{gen-Dirac}
 D\psi \ =\ \sum_{i} Z_i\cdot Z_i(\psi) + Z_i\cdot \tilde{\Lambda}_{\m}(Z_i)
 \psi\,,
 \end{equation}
where  $Z_1,\,\ldots,Z_n$ denotes any orthonormal basis of $\m$.
In the same article, Ikeda states  a criterion for the 
formal self adjointness of this operator. We restate the result  here, since
there is some confusion about the assumptions on the scalar product in the
original version.
\begin{prop}
Let $M=G/H$ be a homogeneous reductive manifold with a homogeneous spin
structure, $\lan \ , \ \ran$ the scalar product on $\m$ induced by the 
Riemannian metric on $M$, and $\nabla$ the $G$ invariant metric connection 
defined 
by some map  $\Lambda_{\m}:\ \m\ra\so(\m)$. Then the Dirac operator $D$
associated with the connection $\nabla$ is formally self adjoint if and
only if for any vector $X\in \m$ and any orthonormal basis 
$Z_1,\,\ldots,Z_n$ of $\m$, one has
 \begin{equation}\tag{$*$} 
 \sum_i \lan \Lambda_{\m}(Z_i)X,Z_i\ran \ =\ \sum_i \lan 
 [Z_i,X]_{\m},Z_i\ran\,.
 \end{equation}
In particular, this condition is always satisfied if the torsion
$T(X,Y,Z)$ is totally skew symmetric. If the metric $\lan \ , \ \ran$ is in
addition naturally reductive, condition $(*)$ is equivalent to
$\sum\Lambda_{\m}(Z_i)Z_i=0$. 
\end{prop}
\begin{proof}
By a result of Friedrich and Sulanke (\cite{Friedrich&S79}), the 
Dirac operator $D^{\nabla}$ associated with any metric connection $\nabla$
is formally self adjoint if and only if the $\nabla$-divergence of any vector
$X$ coincides with its Riemannian divergence,
 \bdm
 \mathrm{div}^{\nabla}(X)\ :=\ \sum_i \lan Z_i, \nabla_{Z_i}X\ran\ =\ 
 \sum_i \lan Z_i,  \nabla^{\mathrm{LC}}_{Z_i}X\ran\ =:\  \mathrm{div}(X) \,,
 \edm
where $\nabla^{\mathrm{LC}}$ denotes the Levi-Civita connection.
But for any vector $X$, the covariant derivatives are related by
 \bdm
 \nabla_{Z_i}X\ =\ \nabla^{\mathrm{LC}}_{Z_i}X + \frac{1}{2}T(Z_i,X)\,,
 \edm
thus equality of divergences holds if and only if
 \bdm
 \sum_{i}\lan T(Z_i,X),Z_i \ran\ =\ 0\,.
 \edm
Inserting the general formula for the torsion and using the fact that 
$\lan\Lambda_{\m}(X)Z_i,Z_i\ran=0$, one checks that this is equivalent
to condition $(*)$. Since $\lan T(Z_i,X),Z_i \ran= T(Z_i,X,Z_i)$,
condition $(*)$ is always fulfilled if the $(0,3)$-tensor $T$ is
totally skew symmetric. Alternatively, one easily deduces equation $(*)$
from the antisymmetry condition in Lemma~\ref{antisymmetry} by a contraction.
Finally, if the metric is naturally reductive, the right-hand side of $(*)$
vanishes, and by the antisymmetry of $\Lambda_{\m}(Z_i)$ one obtains
$\lan X, \sum \Lambda_{\m}(Z_i)Z_i\ran=0$. This finishes the proof.
\end{proof}
\noindent
Returning to the family $\nabla^t$, our aim is
to rewrite the connection term of the Dirac operator in 
equation~(\ref{gen-Dirac}) as an element of the Clifford algebra $\Cl(\m)$.
Basically this amounts  to the identification of $\spin(\m)$ with
the elements of second degree  in $\Cl(\m)$. We implement the Clifford 
relations via $Z_i\cdot Z_j + Z_j\cdot Z_i= -\delta_{ij}$,
in contrast to \cite{Kostant99} (see \cite{Dirac-Buch-00} for notational
details). The following lemma due to Parthasarathy expresses the
lift of the isotropy representation as an element of the Clifford algebra.
\begin{lem}[{\cite[2.1]{Parthasarathy72}}]\label{Parthasarathy}
%
For any element $Y$ in $\h$, one has
 \bdm
 \adtilde(Y)\ =\ \frac{1}{4}\sum_{i,j=1}^n \lan[Y,Z_i],Z_j \ran Z_i\cdot
 Z_j\,. 
 \edm
\qed
\end{lem}
\noindent
Similarly, any skew symmetric map $\Lambda_{\m}(X):\, \m\ra\m$ may be
expanded in the standard basis $E_{ij}$ of $\so(\m)$ as
 \bdm
 \Lambda_{\m}(X)\ =\ \sum_{i<j}\lan \Lambda_{\m}(X)Z_i,Z_j \ran E_{ij}.
 \edm
Since $E_{ij}$ lifts to $Z_i\cdot Z_j/2$ in the Clifford algebra,
we obtain in complete analogy to the Parthasarathy Lemma:
\begin{lem}\label{Parth-lambda-version}
For any map $\Lambda_{\m}:\, \m\ra\so(\m)$, one has
 \bdm
 \tilde{\Lambda}_{\m}(X)\ =\ \frac{1}{2}\sum_{i<j}\lan \Lambda_{\m}(X)
 Z_i,Z_j \ran Z_i\cdot Z_j\ =\ \frac{1}{4}\sum_{i,j}\lan \Lambda_{\m}(X)
 Z_i,Z_j \ran Z_i\cdot Z_j\,.
 \edm
\qed
\end{lem}
\noindent
In particular, the image of $\Lambda^1_{\m}(Z_i)=[Z_i,-]_{\m}$ in $\Cl(\m)$
may be written
 \bdm
 \tilde{\Lambda}^1_{\m}(Z_i)\ =\ \frac{1}{4}\sum_{j,k}\lan[Z_i,Z_j]_{\m},Z_k 
 \ran Z_j\cdot Z_k\,.
 \edm
Thus, by defining the element
 \bdm
 H\ :=\ \sum_{i=1}^n Z_i\cdot \tilde{\Lambda}^1_{\m}(Z_i)\ =\ \frac{1}{4}
 \sum_{i,j,k}\lan[Z_i,Z_j]_{\m},Z_k \ran Z_i\cdot Z_j\cdot Z_k\ =\
 \frac{3}{2}\sum_{i<j<k}\lan[Z_i,Z_j]_{\m},Z_k \ran Z_i\cdot Z_j\cdot Z_k\,,
 \edm
we can rewrite the Dirac operator corresponding to the connection
 $\nabla^t$ from equation~(\ref{gen-Dirac}) as
 \begin{equation}\label{H-Dirac}
 D^t\psi\ =\ \sum_{i} Z_i\cdot Z_i(\psi) + t\cdot H\cdot\psi\,.
 \end{equation}
\begin{NB}\label{H-T}
We identify differential forms with elements of the Clifford
algebra by 
 \bdm
 \sum_{i_1<\ldots<i_r} \omega_{1\ldots r }\, Z_{i_1}\hut\ldots \hut Z_{i_r}
 \lmapsto
 \sum_{i_1<\ldots<i_r} \omega_{1\ldots r }\,Z_{i_1}\cdot\ldots \cdot Z_{i_r}\,.
 \edm
Thus, the torsion form $T^t(X,Y,Z)=(2t-1)\lan [X,Y]_{\m},Z\ran$ induces 
the element 
 \bdm
 T^t\ =\  (2t-1) \sum_{i<j<k}\lan[Z_i,Z_j]_{\m},Z_k 
 \ran Z_i\cdot Z_j\cdot Z_k
 \edm
of the Clifford algebra, which differs from $H$ only by a numerical
factor,
 \bdm
 T^t\ =\ \frac{2(2t-1)}{3}H\,.
 \edm
The simplicity of equation (\ref{H-Dirac}) is the main reason
why we prefer to work with the element $H$ instead of $T^t$.
\end{NB}
%
\subsection{The cubic element $H$, its square and the Casimir operator}
\label{cubic-square-casimir}
%
It is the cubic element $H$ inside the Clifford algebra $\Cl(\m)$ which 
suggested the name "cubic Dirac operator" to B.\ Kostant. We see that
the fact that $H$ is of degree $3$ inside $\Cl(\m)$ does  not depend
on the particular choice for $\Lambda_{\m}$. The square of $H$
will play an eminent role in our considerations, both for a 
Kostant-Parthasarathy type formula and for general vanishing theorems. 
Notice that the square of
\emph{any} element of degree $3$ inside $\Cl(\m)$ has only terms of 
degree zero and $4$.
\begin{prop}\label{H2}
The terms of degree zero and $4$ of $H^2$ are given by
 \bea[*]
(H^2)_0& =& \frac{3}{8}\sum_{i,j} \lan[Z_i,Z_j]_{\m},[Z_i,Z_j]_{\m}\ran\,,\\
(H^2)_4&=&-\frac{9}{2}\sum_{i<j<k<l}\lan Z_i,\Jacm(Z_j,Z_k,Z_l) \ran
Z_i\cdot Z_j\cdot Z_k\cdot Z_l\,. 
 \eea[*]
The first formula is valid for all $n\geq 3$, while the second holds only
for $n\geq 5$. For $n=3,4$, one has $(H^2)_4=0$.
\end{prop}
\begin{proof}
The contributions of degree zero in $H^2$ are exactly the squares
of the summands of $H$. Because of $(Z_i\cdot Z_j\cdot Z_k)^2=1$, we have
 \bdm
(H^2)_0\ =\ \frac{9}{4}\sum_{i<j<k}\lan[Z_i,Z_j]_{\m},Z_k \ran^2\ =\
\frac{9}{24}\sum_{i,j,k}\lan[Z_i,Z_j]_{\m},Z_k \ran\lan[Z_i,Z_j]_{\m},Z_k \ran
\,.
 \edm
For fixed $i,j$, the sum over $k$ is  the coordinate expansion of the
scalar product $\lan[Z_i,Z_j]_{\m},[Z_i,Z_j]_{\m}\ran$, thus
 \bdm
(H^2)_0\ =\ \frac{3}{8}\sum_{i,j}\lan[Z_i,Z_j]_{\m},[Z_i,Z_j]_{\m}\ran\,,
 \edm
as claimed. Contributions of degree $4$ occur if $Z_i\cdot Z_j\cdot Z_k$
is multiplied by $Z_{i'}\cdot Z_{j'}\cdot Z_{k'}$ with exactly one
common index. Since this requires at least $5$ different indices,
it follows that there are no terms of fourth degree for $n\leq 4$. For
the moment, put aside the overall factor $9/4$ of $H^2$.
We explain the occurrence of the term proportional to
$Z_{1234}:=Z_1\cdot Z_2\cdot Z_3\cdot Z_4$ in detail, the others are 
obtained in a 
similar way. Since $H$ consists of ordered tuples proportional to
$Z_{ijk}:=Z_i\cdot Z_j\cdot Z_k$, $i<j<k$, the only way to obtain a term in
$Z_{1234}$ is to multiply $Z_{12k}$ by $Z_{34k}$, $Z_{13k}$ by $Z_{24k}$
and $Z_{14k}$ by $Z_{23k}$ for any index $k\geq 5$. First we notice that 
the order of multiplication is irrelevant, since
 \bdm
 Z_{12k}\cdot Z_{34k}\ =\ Z_{34k}\cdot Z_{12k},\quad
 Z_{13k}\cdot Z_{24k}\ =\ Z_{24k}\cdot Z_{13k},\quad\text{and} \
 Z_{14k}\cdot Z_{23k}\ =\ Z_{23k}\cdot Z_{14k}\,.
 \edm
Every term will thus have multiplicity two. In the next step, these products 
have to be rearranged in order to be proportional to $Z_{1234}$:
 \bdm
 Z_{12k}\cdot Z_{34k}\ =\ - Z_{1234},\quad
 Z_{13k}\cdot Z_{24k}\ =\ +Z_{1234},\quad
 Z_{14k}\cdot Z_{23k}\ =\ - Z_{1234}\,.
 \edm
The total contribution coming from the products   $Z_{12k}$ by $Z_{34k}$
is thus
 \bdm
 (*)\ :=\  -2\, Z_{1234}\sum_{k\geq 5} \lan[Z_1,Z_2]_{\m},Z_k \ran 
 \lan[Z_3,Z_4]_{\m},Z_k \ran\,.
 \edm
This is equal to the sum over all $k$, since the additional terms
are zero. However, it shows that the sum is precisely the
expansion of the scalar product $\lan[Z_1,Z_2]_{\m},[Z_3,Z_4]_{\m}\ran$:
 \bdm
 (*)\ =\ - 2\, Z_{1234}\sum_{k=1}^n \lan[Z_1,Z_2]_{\m},Z_k \ran 
 \lan[Z_3,Z_4]_{\m},Z_k \ran\ =\ - 2\, Z_{1234}
 \lan[Z_1,Z_2]_{\m},[Z_3,Z_4]_{\m}\ran\,.
 \edm
After a similar simplification of the other two contributions, the fourth degree
term in $H^2$ proportional to $Z_{1234}$ is finally equal to
 \bdm
 (**)\ :=\ 2 \left[- \lan[Z_1,Z_2]_{\m},[Z_3,Z_4]_{\m}\ran + 
 \lan[Z_1,Z_3]_{\m}, [Z_2,Z_4]_{\m}\ran - 
 \lan[Z_1,Z_4]_{\m},[Z_2,Z_3]_{\m}\ran  \right]\cdot Z_{1234}\,.
 \edm
This, in turn, may be rewritten as
 \bdm
 (**)\ =\ -2 \lan Z_1, \Jacm(Z_2,Z_3,Z_4)\ran \cdot Z_{1234}\,.
 \edm
Putting back in the factor $9/4$, we get the factor $-9/2$ as stated in
the lemma.
\end{proof}
\noindent
For later reference, we compute the anticommutator of $H$ with an
element $Z_l$ for arbitrary $l$.
\begin{lem}\label{anticomm-H}
For any $l$, one has $\ H\cdot Z_l +Z_l\cdot H\, =\, -\frac{3}{2}
\displaystyle{\sum_{i,j}} \lan Z_l,[Z_i,Z_j]_{\m}\ran Z_i\cdot Z_j$.
\end{lem}
\begin{proof}
By definition,
 \bdm
 H\cdot Z_l+ Z_l\cdot H\, =\,\frac{1}{4} \sum_{i,j,k}
 \lan[Z_i,Z_j]_{\m},Z_k \ran \big(Z_i\cdot Z_j\cdot Z_k\cdot Z_l+Z_l
 \cdot  Z_i\cdot Z_j\cdot Z_k)\,.
 \edm
If all four indices $i,j,k,l$ are pairwise different, 
 \bdm
Z_i\cdot Z_j\cdot Z_k\cdot Z_l\ =\ - Z_l\cdot  Z_i\cdot Z_j\cdot Z_k,
 \edm
and the corresponding summand vanishes. Thus, the sum may be split into those
parts where $l$ is one of the indices $i$, $j$ and $k$,
 respectively:
 \bea[*]
 H\cdot Z_l+ Z_l\cdot H&=& \frac{1}{4} \sum_{j,k}
 \lan[Z_l,Z_j]_{\m},Z_k \ran \big(Z_l\cdot Z_j\cdot Z_k\cdot Z_l+Z_l
 \cdot  Z_l\cdot Z_j\cdot Z_k)\\
 &+& \frac{1}{4} \sum_{i,k}
 \lan[Z_i,Z_l]_{\m},Z_k \ran \big(Z_i\cdot Z_l\cdot Z_k\cdot Z_l+Z_l
 \cdot  Z_i\cdot Z_l\cdot Z_k)\\
 &+& \frac{1}{4} \sum_{i,j}
 \lan[Z_i,Z_j]_{\m},Z_l \ran \big(Z_i\cdot Z_j\cdot Z_l\cdot Z_l+Z_l
 \cdot  Z_i\cdot Z_j\cdot Z_l)\,.
 \eea[*]
We simplify the mixed products to get
 \bea[*]
  H\cdot Z_l+ Z_l\cdot H& =& - \frac{1}{2} \sum_{j,k}
 \lan[Z_l,Z_j]_{\m},Z_k \ran Z_j\cdot Z_k + \frac{1}{2}\sum_{i,k}
 \lan[Z_i,Z_l]_{\m},Z_k \ran Z_i\cdot Z_k\\
 &-& \frac{1}{2}\sum_{i,j}\lan[Z_i,Z_j]_{\m},Z_l \ran Z_i\cdot Z_j\,.
 \eea[*]
Using the invariance property of the scalar product and renaming
the summation indices, this is easily seen to be the desired expression.
\end{proof}
\noindent
Finally, we compute the image of the quadratic Casimir operator of $\h$ inside
the Clifford algebra. Since the $\Ad(G)$ invariant extension $Q$ of
$\lan\ , \ \ran$ is not necessarily positive definite when restricted to $\h$, 
it is more appropriate to work with dual rather than with orthonormal bases.
So  pick bases $X_i,Y_i$ of $\h$ wich are dual with respect to $Q_{\h}$, 
i.\ e., $Q_{\h}(X_i,Y_j)=\delta_{ij}$. The lift of the Casimir operator 
of $\h$ is  defined as
 \bdm
 \Ctilde_{\h}\ =\ - \sum_{i}\adtilde(X_i)\circ \adtilde(Y_i)\,.
 \edm
By the Parthasarathy Lemma (Lemma~\ref{Parthasarathy}), 
 \bdm
 \adtilde(X_i)\ =\ \frac{1}{4}\sum_{j,k} \lan[X_i,Z_j],Z_k \ran Z_j\cdot
 Z_k
 \edm
and similarly for $\adtilde(Y_i)$. Thus,
 \bdm
 \Ctilde_{\h}\ =\ - \frac{1}{16}\sum_{i}\sum_{j,k,l,p}\lan[X_i,Z_j],Z_k \ran 
 \lan[Y_i,Z_l],Z_p \ran Z_j\cdot Z_k\cdot Z_l\cdot Z_p\,.
 \edm
We may get rid of the sum over $i$ immediately. Since $\m$ is
orthogonal to $\h$, we can  rewrite $\Ctilde_{\h}$ as
 \bea[*]
 \Ctilde_{\h}& =& - \frac{1}{16}\sum_{i}\sum_{j,k,l,p} Q([X_i,Z_j],Z_k)
 Q([Y_i,Z_l],Z_p) Z_j\cdot Z_k\cdot Z_l\cdot Z_p\\
 &=& - \frac{1}{16}\sum_{i}\sum_{j,k,l,p} Q(X_i,[Z_j,Z_k])
 Q(Y_i,[Z_l,Z_p]) Z_j\cdot Z_k\cdot Z_l\cdot Z_p\,.
 \eea[*]
For fixed $j,k,l$ and $p$, the sum over $i$ is again the 
expansion of the $\h$ part of $Q([Z_j,Z_k],[Z_l,Z_p])$, yielding
 \begin{equation}\label{Casimir-h-red-form}
 \Ctilde_{\h}\ =\ - \frac{1}{16}\sum_{j,k,l,p}Q_{\h}([Z_j,Z_k],[Z_l,Z_p])
 Z_j\cdot Z_k\cdot Z_l\cdot Z_p\,.
 \end{equation}
This expression has the advantage that it does not contain the dual
bases $X_i,Y_i$ any more.
It turns out that $\Ctilde_{\h}$ has no second degree term, for such a term
would occur if the two index pairs $(j,k)$ and $(l,p)$ had exactly
one common index, for example, $j=l$. But such a term would appear
twice, namely, as $Z_j\cdot Z_k\cdot Z_j\cdot Z_p$ and as
$Z_j\cdot Z_p\cdot Z_j\cdot Z_k$, and these cancel out each other.
\begin{prop}\label{Casimir-h-terms}
The terms of degree zero and $4$ of $\Ctilde_{\h}$ are given for all $n\geq 3$ by
 \bea[*]
(\Ctilde_{\h})_0& =& \frac{1}{8}\sum_{i,j} Q_{\h}([Z_i,Z_j],[Z_i,Z_j])\,,\\
(\Ctilde_{\h})_4&=&-\frac{1}{2}\sum_{i<j<k<l}\lan Z_i,\Jach(Z_j,Z_k,Z_l) \ran
Z_i\cdot Z_j\cdot Z_k\cdot Z_l\,. 
 \eea[*]
In particular, $(\Ctilde_{\h})_4$ vanishes identically for $n\leq 3$, but not
for $n=4$. 
\end{prop}
\begin{proof}
As the form of the result suggests, the proof is  similar to
the computation of $H^2$ (Proposition~\ref{H2}). This is why we 
shall be brief. For the zero degree term, $(j,k)=(l,p)$, and each  term
of this kind appears twice, thus
 \bdm
 (\Ctilde_{\h})_0\ =\ -\frac{1}{8}\sum_{i,j}Q_{\h}([Z_i,Z_j],[Z_i,Z_j])
 Z_i\cdot Z_j\cdot Z_i\cdot Z_j\,.
 \edm
Since $Z_i\cdot Z_j\cdot Z_i\cdot Z_j=-1$, we obtain the first part of the
proposition. For the fourth degree contribution, rewrite the Casimir operator
as
 \begin{equation}
 \Ctilde_{\h}\ =\ - \frac{1}{4}\sum_{j<k,l<p}Q_{\h}([Z_j,Z_k],[Z_l,Z_p])
 Z_j\cdot Z_k\cdot Z_l\cdot Z_p\,.
 \end{equation}
Then the index pairs $(j,k)$ and $(l,p)$ have to be completely disjoint. 
Agein we look  only at the term that is  proportional to 
$Z_{1234}:=Z_1\cdot Z_2\cdot Z_3\cdot Z_4$. It may be obtained
by multiplying $Z_{12}$ by $Z_{34}$, $Z_{13}$ by $Z_{24}$ and $Z_{14}$ by 
$Z_{23}$.
Again, these elements commute, so  we only need to consider each
product in the order of multiplication just given  and count it twice.
Restoring the order of indices in these products, one sees that
the term in $(\Ctilde_{\h})_4$ proportional to $Z_{1234}$ looks like
 \bdm
 (*)\ :=\ - \frac{2}{4} \left[Q_{\h}([Z_1,Z_2],[Z_3,Z_4])-
 Q_{\h}([Z_1,Z_3],[Z_2,Z_4]) + Q_{\h}([Z_1,Z_4],[Z_2,Z_3]) 
 \right]\cdot Z_{1234}\,.
 \edm
By the properties of $Q$, the first scalar product may be formulated
differently:
 \bdm
 Q_{\h}([Z_1,Z_2],[Z_3,Z_4])\ =\ Q([Z_1,Z_2],[Z_3,Z_4]_{\h})\ =\
 Q(Z_1,[Z_2,[Z_3,Z_4]_{\h}])\,.
 \edm
Rewriting the other two products in a similar way, we see that
 \bdm
 (*)\ =\  - \frac{1}{2}\,Q(Z_1,\Jach(Z_2,Z_3,Z_4) )\cdot Z_{1234} \,.
 \edm
\end{proof}
%
\subsection{A Kostant-Parthasarathy type formula for $(D^t)^2$}
%
If $M=G/H$ is a symmetric space, it is well known that besides the
general Schr\"odinger-Lichnerowicz formula for $D^2$, which is valid on any
Riemannian manifold, there exists a formula expressing $D^2$ in
terms of Casimir operators due to Parthasarathy (see also 
\cite[Remark 1.63]{Kostant99}). The Dirac operator $D$ is defined relative 
to the Levi-Civita connection, which coincides with our one-parameter family
$\nabla^t$, and $\lan\ , \ \ran$ denotes an $\Ad(G)$ invariant scalar
product on $\g$ whose restriction to $\m$ is positive definite.
Let $\Scal$  be the scalar curvature of the symmetric space $M$ and  
$\Omega_{G}$ the Casimir operator of $G$, viewed as a second order 
differential operator. 
\begin{thm}[{\cite[Prop.3.1]{Parthasarathy72}}, {\cite[Ch. 3]{Dirac-Buch-00}}]
\label{Kos-Parth-D2-symm}
On a symmetric space $M=G/H$, one has 
\bdm
 D^2\ =\ \Omega_{G} + \frac{1}{8}\Scal\,,
\edm
and the scalar curvature may be rewritten as
 $\Scal\, =\, 8\cdot (\lan\vrho_{\g},\vrho_{\g}\ran -\lan\vrho_{\h},\vrho_{\h}
 \ran )$.
\end{thm}
\noindent
This formula is the starting point for vanishing theorems, the
realization of discrete series representations in the kernel of $D$,
and it allows the computation of the full spectrum of $D$ on $M$.
If we now go back to the situation studied in this article, i.\,e.,
a reductive homogeneous space $G/H$ endowed with a naturally reductive
metric $\lan\ , \ \ran$ on $\m$, then, a priori, the steps in the proof
of Theorem~\ref{Kos-Parth-D2-symm} cannot be performed any longer.
To prove a Kostant-Parthasarathy type formula in this situation, we
recall the general expression for the Dirac operator associated with
the connection $\nabla^t$ from equation~(\ref{H-Dirac}) and split it into 
the terms coming from the canonical connection and the $3$-form
$H$, respectively:
 \begin{equation}
 D^t\psi \ =\ \sum_{i} Z_i\cdot Z_i(\psi) + Z_i\cdot 
 \tilde{\Lambda}^t_{\m}(Z_i) \psi\ =:\ D^0\psi+ D^t_H\psi\,.
 \end{equation}
First, notice that the equivariance property of spinors  implies that 
the action on spinors of vector fields coming from $\m$ is by ``true'' 
differential operators, while the action of vector fields in $\h$ is in 
fact purely algebraic.
\begin{lem}\label{H-on-spinors}
Let $\psi$ be a spinor, i.\,e., a section in $S=G\x_{\Adtilde}\Delta_{\m}$ and 
$X$ an element of $\h$, identified with the left invariant vector field it 
induces. Then
 \bdm
 X(\psi) \ =\ -\adtilde(X)\cdot \psi\,,
 \edm
where $\adtilde(X)\cdot \psi$ denotes the Clifford product of the spinor $\psi$
with the element $\adtilde(X)\subset\spin(\m)\subset\Cl(\m)$.
\end{lem}
\begin{proof}
We identify $\psi$ with a map  $\psi:\, G\ra\Delta_{\m}$ such that 
$\psi(gh)\ =\ \kappa(\Adtilde(h^{-1}))\psi(g)$ for all $g\in G$ and $h\in H$. 
Then one has
 \bdm
 X\psi(g)\ =\ \frac{d}{ds} \psi(ge^{sX})\big|_{s=0}\ =\
 \frac{d}{ds}\kappa(\Adtilde(e^{-sX}))\psi(g)\big|_{s=0}\ =\ -\kappa
 (\adtilde(X))\psi(g)\,.
 \edm
Thus, $X(\psi)=-\kappa(\adtilde(X))\psi=-\adtilde(X)\cdot \psi$, as claimed.
\end{proof}
\begin{NB}
In \cite[Section 2]{Kostant99} and \cite[Chapter 10.5]{Sternberg99}, the map 
assigning to $X\in \h$ the sum
 \bdm
 X(-) + \adtilde(X)\cdot - 
 \edm
is called the ``diagonal'' map from $\h$ to $\mathcal{U}(\g)\ox\Cl(\m)$.  
The assumption that the action is on spinors thus implies that this
diagonal map is equal to zero. In particular, the diagonal Casimir
operator of $\h$ vanishes in the formula for $(D^t)^2$.
\end{NB}
\begin{prop}\label{D0-square}
The square of $D^0$, the Dirac operator corresponding to the
canonical connection,  is given by
 \bdm
 (D^0)^2\psi\ =\  -\sum_{i}Z_i^2(\psi)+ 2\, \Ctilde_{\h} + 
 \frac{1}{2}\sum_{i,j,k}\lan [Z_i,Z_j]_{\m},Z_k\ran Z_i\cdot Z_j
 \cdot Z_k(\psi) \,.
 \edm
\end{prop}
\noindent
Before proceeding to the proof, let us make a short remark on how this formula
is to be understood. In the first term, one has to take  the
derivative of $\psi$ along all vector fields $Z_i$ twice, thus yielding a
second order differential operator. By $\Ctilde_{\h}$, we mean the image of the
Casimir operator of $\h$ inside $\Cl(\m)$ as described in Section
\ref{cubic-square-casimir}. Finally, $Z_i\cdot Z_j\cdot$ denotes the 
Clifford product of $Z_i$ and $Z_j$, whereas $Z_k$ acts again as a 
derivative. Thus the last term is a first order differential operator. Notice 
that Clifford multiplication by any constant element in $\Cl(\m)$
commutes with derivation along $\m$.
\begin{proof}
We compute $(D^0)^2$ as follows:
 \bdm
 (D^0)^2\psi\ =\ \sum_iZ_i\cdot Z_i(\sum_j Z_j\cdot Z_j(\psi))\ =\
 \sum_{i,j}Z_i\cdot Z_j\cdot (Z_i Z_j (\psi))\,.
 \edm
We divide the sum into the diagonal ($i=j$) and off-diagonal ($i\neq j$)
terms and see that this separates the second and the first order differential
operator contribution,
 \bdm
 (D^0)^2\psi\ =\ -\sum_i Z_i^2(\psi) + \frac{1}{2} \sum_{i,j}
 Z_i\cdot Z_j\cdot [Z_i,Z_j](\psi)\,.
 \edm
We concentrate our attention on the second term. Split the commutator
into its $\m$ and $\h$ part, then write the $\m$ part again in the
orthonormal basis $Z_1,\ldots,Z_n$ to obtain
 \bea[*]
 \frac{1}{2} \sum_{i,j}
 Z_i\cdot Z_j\cdot [Z_i,Z_j](\psi)&=& \frac{1}{2} \sum_{i,j}
 Z_i\cdot Z_j\cdot ([Z_i,Z_j]_{\m}(\psi)+[Z_i,Z_j]_{\h}(\psi))\\
 &=& \frac{1}{2} \sum_{i,j,k}\lan Z_k, [Z_i,Z_j]_{\m} \ran 
 Z_i\cdot Z_j\cdot Z_k(\psi) + \frac{1}{2} \sum_{i,j}
 Z_i\cdot Z_j\cdot [Z_i,Z_j]_{\h}(\psi)\,.
 \eea[*]
This takes care of the last term in the formula of Proposition~\ref{D0-square}.
Thus it remains to show that
 \bdm
 \frac{1}{2} \sum_{i,j}
 Z_i\cdot Z_j\cdot [Z_i,Z_j]_{\h}(\psi)\ =\ 2\, \Ctilde_{\h}\,.
 \edm
The action of the commutators $[Z_i,Z_j]_{\h}$ on the spinor $\psi$ is
first transformed into Clifford multiplication  by the adjoint representation
as explained in Lemma~\ref{H-on-spinors}, then rewritten
in terms of an orthonormal basis according to the Parthasarathy Lemma 
(Lemma~\ref{Parthasarathy}),
 \bea[*]
 \frac{1}{2} \sum_{i,j} Z_i\cdot Z_j\cdot [Z_i,Z_j]_{\h}(\psi)&=&
 - \frac{1}{2} \sum_{i,j} Z_i\cdot Z_j\cdot\adtilde([Z_i,Z_j]_{\h})\cdot\psi\\
 & =& - \frac{1}{8} \sum_{i,j}Z_i\cdot Z_j\sum_{p,q}
 \lan [[Z_i,Z_j]_{\h},Z_p],Z_q\ran Z_p\cdot Z_q\cdot \psi\,.
 \eea[*]
But since $\lan [[Z_i,Z_j]_{\h},Z_p],Z_q\ran=Q_{\h}([Z_i,Z_j],[Z_p,Z_q])$,
this is $2\, \Ctilde_{\h}$ by equation~(\ref{Casimir-h-red-form}).
\end{proof}
\noindent
With the preparations of Section~\ref{cubic-square-casimir}, the other two
terms in the expression for $(D^t)^2$ are relatively easy to compute.
We  denote the Casimir operator of the full Lie algebra $\g$ by
$\Omega_{\g}$,
 \bdm
 \Omega_{\g}\psi \ =\ - \sum_{i} Z^2_i(\psi) + \Ctilde_{\h}\cdot\psi\,.
 \edm
We decided to use a symbol different from $C$ in order
to emphasize that $\Omega_{\g}$ is a real second order differential
operator, as opposed to  $\Ctilde_{\h}$, which is a constant element of the
Clifford algebra. In particular, the result of Lemma~\ref{D0-square}
may be restated as
 \begin{equation}\label{D0-square-bis}
 (D^0)^2\psi\ =\ \Omega_{\g} + \Ctilde_{\h} + 
 \frac{1}{2}\sum_{i,j,k}\lan [Z_i,Z_j]_{\m},Z_k\ran Z_i\cdot Z_j
 \cdot Z_k(\psi) \,.
 \end{equation}
First we state the formula in its most general form.
\begin{thm}[General Kostant-Parthasarathy formula]\label{K-P-1}
For $n\geq 5$, the square of $D^t$ is given by
 \bea[*]
 (D^t)^2\psi & =& \Omega_{\g}(\psi) + \frac{1}{2}(1-3t) \sum_{i,j,k}
 \lan [Z_i,Z_j]_{\m},Z_k \ran Z_i\cdot Z_j\cdot Z_k(\psi)\\
 &-& \frac{1}{2}\sum_{i<j<k<l}\lan Z_i, \Jach(Z_j,Z_k,Z_l)+ 9t^2
 \Jacm(Z_j,Z_k,Z_l)\ran\cdot Z_i\cdot Z_j\cdot Z_k\cdot Z_l\cdot\psi \\
 &+&  \frac{1}{8}\sum_{i,j} Q_{\h}([Z_i,Z_j],[Z_i,Z_j])\psi
 +\frac{3}{8}t^2 \sum_{i,j} Q_{\m}([Z_i,Z_j],[Z_i,Z_j])\psi\,.
 \eea[*]
For $n\leq 4$, one has 
 \bea[*]
 (D^t)^2\psi & =& \Omega_{\g}(\psi) + \frac{1}{2}(1-3t) \sum_{i,j,k}
 \lan [Z_i,Z_j]_{\m},Z_k \ran Z_i\cdot Z_j\cdot Z_k(\psi)\\
 &-& \frac{1}{2}\sum_{i<j<k<l}\lan Z_i, \Jach(Z_j,Z_k,Z_l)
 \ran\cdot Z_i\cdot Z_j\cdot Z_k\cdot Z_l\cdot\psi \\
 &+&  \frac{1}{8}\sum_{i,j} Q_{\h}([Z_i,Z_j],[Z_i,Z_j])\psi
 +\frac{3}{8}t^2 \sum_{i,j} Q_{\m}([Z_i,Z_j],[Z_i,Z_j])\psi\,.
 \eea[*]
\end{thm}
\begin{proof}
The mixed term is the first order differential operator
 \bea[*]
 (D^0 D^t_H+D^t_HD^0)\psi & =& t \sum_p\left[Z_p\cdot Z_p(H\cdot\psi)
 +H\cdot Z_p\cdot Z_p(\psi) \right]\\
 &=& t \sum_p\left[Z_p\cdot H +H\cdot Z_p\right]\cdot Z_p(\psi)\,.
 \eea[*]
In Lemma~\ref{anticomm-H}, we computed the anticommutator of $H$ with
the vector $Z_p$, which leads us to 
 \bdm
 (D^0 D^t_H+D^t_HD^0)\psi\ =\ - \frac{3}{2}t \sum_p 
 \bigg[\sum_{i,j}\lan Z_p,[Z_i,Z_j]_{\m}\ran Z_i\cdot Z_j\bigg] Z_p(\psi)\,.
 \edm
By Lemma~\ref{H2}, we have
 \bdm
(D^t_H)^2\psi\ =\ -\frac{9}{2}t^2 \sum_{i<j<k<l}\lan Z_i,\Jacm(Z_j,Z_k,Z_l) 
\ran Z_i\cdot Z_j\cdot Z_k\cdot Z_l\cdot\psi
+ \frac{3}{8}t^2\sum_{i,j} \lan[Z_i,Z_j]_{\m},[Z_i,Z_j]_{\m}\ran\psi
 \edm
for $n\geq 5$ and
 \bdm
(D^t_H)^2\psi\ =\ \frac{3}{8}t^2\sum_{i,j} \lan[Z_i,Z_j]_{\m},[Z_i,Z_j]_{\m}
\ran\psi
 \edm
otherwise. Together with equation~(\ref{D0-square-bis}) and the formula 
for $\Ctilde_{\h}$ from Proposition~\ref{Casimir-h-terms}, one obtains 
the desired formulas.
\end{proof}
\noindent
Now it becomes clear that the particular choice $t=1/3$ leads to substantial
simplifications in case of $n=3$ or $n\geq 5$. In fact, the second part of the 
first line vanishes identically, the second line is zero by the Jacobi 
identity in $\g$ ($n\geq 5$) or for dimensional reason ($n=3$), and the scalar 
contributions in the last line appear in a very precise
ratio, which will allow some further simplification. It is a strange effect
that no  simplification is possible for $n=4$.
\begin{thm}[The Kostant-Parthasarathy formula for $t=1/3$]\label{K-P-2}
For $n=3$ or $n\geq 5$ and $t=1/3$, the general formula for $(D^t)^2$ 
reduces to
 \bea[*]
 (D^{1/3})^2\psi & = & \Omega_{\g}(\psi)+ \frac{1}{8} \bigg[ 
 \sum_{i,j} Q_{\h}([Z_i,Z_j],[Z_i,Z_j]) + 
 \frac{1}{3}\sum_{i,j} Q_{\m}([Z_i,Z_j],[Z_i,Z_j]) \bigg]\psi \\
 & =& \Omega_{\g}(\psi) + \frac{1}{8}
 \bigg[\Scal^{1/3}+ \frac{1}{9}\sum_{i,j}Q_{\m}([Z_i,Z_j],[Z_i,Z_j])\bigg]
 \psi\,.
 \eea[*]
\qed
\end{thm}
\begin{NB}
In particular, one  immediately recovers the classical Parthasarathy formula
for a symmetric space (Theorem~\ref{Kos-Parth-D2-symm}), since then 
all scalar curvatures coincide and $[Z_i,Z_j]\in \h$.
\end{NB}
\noindent
As in the classical Parthasarathy formula, the scalar term as well as the
eigenvalues of $\Omega_{\g}(\psi)$ may be expressed
in representation theoretical terms if $G$  (and hence $M$) is compact. 
Consider the unique $\Ad(G)$ invariant
extension $Q$ of the scalar product $\lan\,,\,\ran$ on $\m$ to the full
Lie algebra $\g$, which exists by Kostant's Theorem. Thus, $Q$ is a 
multiple of  the Killing form on any simple factor of $\g$; however, $Q$  is 
not necessarily positive definite, hence the scaling factors may be of
different sign. If they are such that  $Q$ is positive definite,
the metric $\lan\,,\,\ran$ is said to be \emph{normal homogeneous}.\\

\noindent
We begin with a more careful analysis of the Casimir operator 
$\Omega_{\g}(\psi)$. By the same arguments as in the symmetric space
case, $\Omega_{\g}(\psi)$ is a $G$ invariant differential operator,
and this property does not depend on the signs of $Q$. We
sketch the argument for completeness: On every simple
summand $\g_i$ of $\g$, $Q_i:=Q\vert_{\g_i}$ is either a positive or a negative
multiple of the Killing form, and $\Ad(g)$ maps $\g_i$ into itself.
Hence, in either case, the adjoint action of $G$ transforms an orthonormal 
base $\tilde{Z}_1,\ldots,\tilde{Z}_m$
of $\g_i$ into an orthonormal base, and dual bases $\tilde{X}_1,\tilde{Y}_1,
\ldots,\tilde{X}_m,\tilde{Y}_m$ of $\g_i$ are mapped 
to dual bases: 
 \bdm
 Q_i(\Ad(g)\tilde{Z}_k,\Ad(g)\tilde{Z}_l)\ =\ Q_i(\tilde{Z}_k,\tilde{Z}_l), 
 \quad
 Q_i(\Ad(g)\tilde{X}_k,\Ad(g)\tilde{Y}_l)\ =\ Q_i(\tilde{X}_k,\tilde{Y}_l)\,.
 \edm
Now consider the Frobenius decomposition of the square integrable spinors
into irreducible finite-dimensional representations $V_{\lambda}$ of $G$,
 \bdm
 L^2(S)\ =\ \sum_{\lambda\in\hat{G}} M_{\lambda}\ox V_{\lambda},
 \edm
where $M_{\lambda}$ denotes the multiplicity space of $V_{\lambda}$.
Let $\vrho_{\lambda}:\, G\ra \GL(V_{\lambda}) $ be the representation
with highest weight $\lambda$, and $d\vrho_{\lambda}$ its differential.
Then $\Omega_{\g_i}$ acts on $V_{\lambda}$ by
\bdm
d\vrho_{\lambda}(\Omega_{\g_i})\ =\ - \sum_{k=1}^m d\vrho_{\lambda}
(\tilde{Z}_k)^2 \ \ \text{ or }\ \ d\vrho_{\lambda}(\Omega_{\g_i})\ =\ - 
\sum_{k=1}^m d\vrho_{\lambda}(\tilde{X}_k)d\vrho_{\lambda}(\tilde{Y}_k)\,.
\edm
However, for any element $X\in\g_i$, one checks immediately
 \bdm
 \vrho_{\lambda}(g) d\vrho_{\lambda}(X)\vrho_{\lambda}(g^{-1})\ =\
 d\vrho_{\lambda}(\Ad(g)X),
 \edm
hence $\Omega_{\g_i}$ commutes with the action of $g\in G$ on 
$V_{\lambda}$, as claimed. Furthermore, it acts by multiplication
by the well-known  eigenvalue
 \bdm
 Q_i(\lambda+\vrho_i,\lambda+\vrho_i ) - Q_i(\vrho_i,\vrho_i ),
 \edm
whose sign, however, depends on whether $Q_i$ is a positive or a negative
multiple of the Killing form on $\g_i$. Here, $\vrho_i$ denotes the
half sum of positive roots of $\g_i$, as usually. 
Since the center of $G$ does
not contribute to the total eigenvalue of $\Omega_{\g}$, we conclude:
\begin{lem}\label{non-negative}
The operator $\Omega_{\g}$ is non negative if
the metric $\lan\,,\,\ran$ is normal homogeneous or if the
negative definite contribution to $Q$ comes from an abelian summand
in $\g$.\qed
\end{lem}
\noindent
In a forthcoming paper, we will discuss examples
where $Q$ has also a simple summand on which $Q$ is negative definite
and show that $\Omega_{\g}$ has  negative eigenvalues.
We use these remarks to express the scalar term in Theorem \ref{K-P-2} in
a different way.
\begin{lem}\label{comp-scalar-in-K-P}
Let $G$ be compact, $n=3$ or $n\geq 5$, and denote by $\vrho_{\g}$
and $\vrho_{\h}$ the half sum of the positive roots of $\g$ and $\h$, 
respectively. Then the Kostant-Parthasarathy
formula for $(D^{1/3})^2$ may be restated as
 \bdm
 (D^{1/3})^2\psi \ = \ \Omega_{\g}(\psi)+ \left[Q( \vrho_{\g},\vrho_{\g})-
 Q(\vrho_{\h},\vrho_{\h})\right]\psi \ =\ \Omega_{\g}(\psi)+
 \lan \vrho_{\g}- \vrho_{\h},\vrho_{\g}- \vrho_{\h}  \ran
 \psi\,.
 \edm
In particular, the scalar term is positive independently
of the properties of $Q$.
\end{lem}
\begin{proof}
Consider the eightfold multiple of the term under consideration
and regroup it as
 \bea[*]
 8 ((D^{1/3})^2-\Omega_{\g})& = & 
 \sum_{i,j} Q_{\h}([Z_i,Z_j],[Z_i,Z_j]) + 
 \frac{1}{3}\sum_{i,j} Q_{\m}([Z_i,Z_j],[Z_i,Z_j]) \\
 &=&  \frac{1}{3}\bigg[\sum_{i,j} Q([Z_i,Z_j],[Z_i,Z_j])
  + 2  \sum_{i,j} Q_{\h}([Z_i,Z_j],[Z_i,Z_j])\bigg]\,.
 \eea[*]
The first summand can easily be seen to be a trace over $\m$,
\bdm
 \sum_{i,j} Q([Z_i,Z_j],[Z_i,Z_j])\, =\, - \sum_{i,j} Q([Z_i,[Z_i,Z_j]],Z_j)
\, =\, -\sum_j Q( \sum_i (\ad Z_i)^2,Z_j)\, =\, -\tr_{\m}\sum_i (\ad Z_i)^2\,.
\edm
For the second term, we first notice that it may be rewritten by expanding
and contracting in two different ways as
\bea[*]
\sum_{i,j} Q_{\h}([Z_i,Z_j],[Z_i,Z_j])& =& 
\sum_{i,j,k} Q(X_k,[Z_i,Z_j])Q(Y_k,[Z_i,Z_j])\, =\,
\sum_{i,j,k} Q([X_k,Z_i],Z_j)Q([Y_k,Z_i],Z_j)\\ 
&=& \sum_{i,k} Q([X_k,Z_i],[Y_k,Z_i]). 
\eea[*]
This, in turn, can be identified with two different kinds of traces:
On the one hand, this is
\bdm
- \sum_{i,k} Q([Z_i,[Z_i,X_k]],Y_k)\ =\ - \tr_{\h}\sum_i (\ad Z_i)^2\,,
\edm
on the other hand, this reads
 \bdm
- \sum_{i,k} Q([X_k,[Y_k,Z_i]],Z_i)\ =\ -\tr_{\m}\sum_k (\ad X_k)(\ad Y_k)
\ =\  \tr_{\m}C_{\h}\,,
 \edm
were $C_{\h}$ denotes the ``unlifted'' Casimir operator of $\h$, i.\,e., 
its usual action on $\g$ via the adjoint
representation. Now, since the sum we have just treated appears twice, we use 
each way of writing it once to obtain
 \bea[*]
 8 ((D^{1/3})^2-\Omega_{\g})& = & \frac{1}{3}\bigg[
-\tr_{\m}\sum_i (\ad Z_i)^2 - \tr_{\h}\sum_i (\ad Z_i)^2 + \tr_{\m}C_{\h}
\bigg]\\
& =&  \frac{1}{3}\bigg[-\tr_{\g}\sum_i (\ad Z_i)^2 + \tr_{\g}C_{\h}
- \tr_{\h}C_{\h}\bigg]\\
&=& \frac{1}{3}\bigg[ \tr_{\g}C_{\g} - \tr_{\h}C_{\h}\bigg]\,.
\eea[*]
Again, $C_{\g}$ is not to be confused with the action of the Casimir
operator of $\g$ on spinors. By looking separately on every simple summand 
where $Q$ is just a multiple of the Killing form, one easily sees
that these traces are the rescaled lengths of the half sum of positive 
roots,
\bdm
\tr_{\g}C_{\g}\ =\ 24\, Q(\vrho_{\g},\vrho_{\g}),
\edm
and similarly for $\h$ (Proposition $1.84$ in \cite{Kostant99}). This
proves the formula. To see that the scalar 
 is positive even for non normal homogeneous metrics,
decompose $\vrho_{\g}= \vrho_{\h}+R$, where $R\in\m$. Since $\m$ and
$\h$ are orthogonal with respect to $Q$, one obtains
 \bdm
 Q(\vrho_{\g},\vrho_{\g}) - Q(\vrho_{\h},\vrho_{\h})\ =\
 Q(\vrho_{\h}+R,\vrho_{\h}+R)  - Q(\vrho_{\h},\vrho_{\h})\ =\
 Q(R,R)\ =\ \lan R,R\ran > 0,
 \edm
since by dimensional reasons $R\neq 0$  and the
scalar product on $\m$ is positive definite.
\end{proof}
\noindent
We can formulate our first conclusion from Theorem \ref{K-P-2}:
\begin{cor}\label{eigenvalue-estimate}
If the operator $\Omega_{\g}$ is non negative,  the first eigenvalue
$\lambda_1^{1/3}$ of the Dirac operator $D^{1/3}$ satisfies the inequality
 \bdm
 \big(\lambda_1^{1/3}\big)^2\ \geq \ Q( \vrho_{\g},\vrho_{\g})-
 Q(\vrho_{\h},\vrho_{\h}) \,.
 \edm
Equality occurs if and only if there exists an algebraic spinor in
$\Delta_{\m}$  which is fixed under the lift $\kappa(\Adtilde H)$ of the isotropy
representation.
\end{cor}
\begin{proof}
By our assumption on $\Omega_{\g}$, its eigenvalue on a spinor $\psi$
can be zero if and only if the Casimir eigenvalue of every simple summand
$\g_i$ of $\g$ vanishes, hence $\psi$ has to lie in the trivial
$G$-representation and is thus constant.
\end{proof}
\noindent
We shall discuss examples of equality at the end of  Section~\ref{examples}.

\begin{NB}\label{D-new-op}
Since $D^t$ is a $G$-invariant differential operator on $M$ by construction,
Theorem \ref{K-P-1} implies that the linear combination
of the first order differential operator and the multiplication
by the element of degree four in the Clifford algebra
appearing in the formula for $(D^t)^2$ is again $G$ invariant for all
$t$. Hence, the first order differential operator
\bdm
\mathcal{D}\psi \ :=\ \sum_{i,j,k} \lan [Z_i,Z_j]_{\m},Z_k \ran 
Z_i\cdot Z_j\cdot Z_k(\psi)
\edm
has to be a  $G$ invariant differential operator, a fact that cannot be seen
directly by any simple arguments. It has no analogue on symmetric spaces and
certainly deserves  further separate investigations.
\end{NB}
%

\section{The equations of type II string theory on naturally reductive spaces}
\subsection{The field equations}
The common sector of type II string theories may be geometrically
described as a tuple $(M^n, \lan\,,\,\ran, H,\Phi, \Psi)$ consisting of a
manifold $M^n$ with a Riemannian metric $\lan\,,\,\ran$, a $3$-form $H$,
a so-called dilaton function $\Phi$ and a spinor field $\Psi$ satisfying
the coupled system of field equations
 \bdm
 \Ric^{\LC}_{ij} - \frac{1}{4}H_{imn}H_{jmn}+ 2\nabla^{\LC}_{i}\del_j\Phi
 \, =\, 0,
 \, \delta(e^{-2\Phi}H)\, =\, 0, \,
 (\nabla^{\LC}_X+ \frac{1}{4}X\haken H) \Psi \, =\, 0,\,
 (d\Phi - \frac{1}{2}H) \Psi \, =\, 0.
 \edm
The first equation generalizes the Einstein equation, the second is
a conservation law, while the first of the spinorial field equations
suggests that the $3$-form $H$  should be the torsion of some metric connection
$\nabla$  with totally skew-symmetric torsion tensor $T=H$. Then the equations 
may be rewritten in terms of $\nabla$:
\bdm
\Ric^{\nabla}+\frac{1}{2}\delta(T)+2\nabla^{\LC} d\Phi\ =\ 0,\quad
\delta(T)\ =\ 2\cdot d\Phi^{\#}\haken T,\quad \nabla\Psi\ =\ 0,\quad
(d\Phi - \frac{1}{2}T)\cdot \Psi \ =\ 0\,.
\edm
If the dilaton  $\Phi$ is constant, the equations may be
simplified even further,
 \bdm
 \Ric^{\nabla}\ =\ 0,\quad \delta(T)\ =\ 0,\quad \nabla\Psi\ =\ 0,\quad
 T\cdot \Psi \ =\ 0\,.
 \edm
In particular, the last equation becomes a purely algebraic condition.
The number of preserved supersymmetries depends essentially on the number
of $\nabla$-parallel spinors. 
For a general background on these equations, we refer to the
article by A.\ Strominger where they appeared first \cite{Strominger86}.
By Lemma~\ref{dT}, we conclude that the second equation is always satisfied
for the family of connections $\nabla^t$.\\

\noindent
Before proceeding further, we add a general observation which
follows easily from the formulas in \cite{Friedrich&I01} and which was
pointed out to us by Bogdan Alexandrov.
\begin{thm}
Let $M^n$ be a compact Riemannian manifold with metric $\lan\,,\,\ran$
and a metric connection $\nabla$  with totally skew symmetric torsion $T$.
Suppose that there exists a spinor field $\psi$ such that all the equations
\bdm
 \Ric^{\nabla}\ =\ 0,\quad \delta(T)\ =\ 0,\quad \nabla\Psi\ =\ 0,\quad
 T\cdot \Psi \ =\ 0
\edm
hold. Then $T=0$ and $\nabla$ is the Levi-Civita connection.
\end{thm}
\begin{proof}
If $\psi$ is $\nabla$-parallel, the Riemannian Dirac operator $D^{\LC}$
acts on  $\psi$ by $D^{\LC}\psi= - 3 T\cdot\psi/4$. The last equation
thus implies $D^{\LC}\psi=0$. By the classical Schr\"odinger-Lichnerowicz
formula, 
\bdm
0\  =\  \int_{M^n} ||\nabla^{\LC} \psi||^2 dM^n + \frac{1}{4}
\int_{M^n}\Scal^{\LC}||\psi||^2 dM^n\,.
\edm
On the other hand, the two Ricci tensors are related by the equation
 \bdm
 \Ric^{\LC}(X,Y)\ =\ \Ric^{\nabla}(X,Y)+\frac{1}{2}(\delta T)(X,Y)+\frac{1}{4}
 \sum_{i=1}^n \lan T(X,e_i),T(Y,e_i)\ran\,,
 \edm
where $e_1,\,\ldots,e_n$ denotes an orthonormal basis. If $\Ric^{\nabla}=0$
and $\delta T=0$, this implies that the Riemannian scalar curvature
is non negative and given by
\bdm
4\, \Scal^{\LC}\ =\ \sum_{i,j=1}^n \lan T(e_i,e_j),T(e_i,e_j)\ran\,.
\edm
Consequently, the scalar curvature $\Scal^{\LC}$ has to vanish
identically, and the torsion form $T$ is zero, too.
\end{proof}
\noindent
Hence, compact solutions to all equations have to be Calabi-Yau
manifolds in dimensions $4$ and $6$, Joyce manifolds in dimensions
$7$ and $8$ etc.

\subsection{Some particular spinor fields}
Consider the situation that the lift of the isotropy representation 
$\kappa(\Adtilde H)$ contains the trivial representation, i.\,e., an
algebraic spinor $\psi$ that is fixed under the action of $H$.
Any such spinor induces a section of the spinor bundle 
$S=G\x_{\kappa(\Adtilde)}\Delta_{\m}$ if viewed as a constant map $G\ra \Delta_{\m}$
and is thus of particular interest.
\begin{thm}\label{constant-fields}
\begin{enumerate}
\item[]
\item Any constant spinor field $\psi$ satisfies the equation
\bdm
\nabla^t_Z \psi\ =\ \frac{t}{3} (Z\haken H)\psi\,.
\edm
In particular, it is parallel with respect to the canonical connection $(t=0)$.
Conversely, any spinor field $\psi$ satisfying $\nabla^0\psi=0$ is
necessarily constant.
\item Any constant spinor field $\psi$ is an eigenspinor of the square of the
Dirac operator $(D^t)^2$, and its eigenvalue does not depend of the
special choice of $\psi$:
 \bdm
 (D^t)^2\psi\ =\ 9t^2 \big[Q(\vrho_{\g},\vrho_{\g}) -
 Q(\vrho_{\h},\vrho_{\h}) \big]\,\psi\,.
 \edm
In particular, $H\cdot\psi\neq 0$ and hence the last string equation can
never hold for a constant spinor.
\end{enumerate}
\end{thm}
\begin{proof}
For a constant spinor field, the formula for the covariant derivative of a 
spinor field (equation \ref{gen-cov-der}) reduces to
$ \nabla^t_Z\psi = 0+ \tilde{\Lambda}^t_{\m}(Z)\psi$.
By Lemma~\ref{Parth-lambda-version},  $\tilde{\Lambda}^t(Z)$ may be
expressed in terms of an orthonormal basis as 
 \bdm
 \nabla^t_Z\psi\ =\ \frac{t}{2}\sum_{j<k}\lan[Z,Z_j]_{\m},Z_k \ran
 Z_j\cdot Z_k\cdot \psi\,.
 \edm
By the definition of $H$, this is easily seen to be $t(Z\haken H)/3$.
Conversely, assume that $\psi$ is parallel
with respect to the canonical connection, i.\,e.\, $Z_i(\psi)=0$ for all $i$.
Then $[Z_i,Z_j](\psi)=0$, and the commutator $[Z_i,Z_j]$ may be split into
its $\m$ and $\h$ part. But the $\m$ part acts again
trivially on $\psi$, hence we obtain
 \bdm
 [Z_i,Z_j]_{\h}(\psi)\ =\ 0\,.
 \edm
By Assumption \ref{assume-transitivity}, $[\m,\m]$ spans all of $\h$,
hence $\h$ also acts trivially on $\psi$, which finishes the argument.
For the second part of the Theorem, we use that
the Dirac operator on a constant spinor is given by 
$D^t\psi=t\,H\cdot\psi$ for any $t$. Since any constant spinor lies 
in the trivial $G$-representation in the Frobenius decomposition of 
$\Gamma(S)$, the eigenvalue of $\Omega_{\g}$ on $\psi$ is zero.
For $t=1/3$, the Kostant-Parthasarathy formula (Theorem~\ref{K-P-2}) 
thus yields
 \bdm
 (D^{1/3})^2\psi\ =\ \big[Q( \vrho_{\g},\vrho_{\g}) - 
 Q( \vrho_{\h},\vrho_{\h})\big] \psi \ =\ \frac{1}{9}H^2\psi\,.
 \edm
This may be understood as a formula for $H^2\psi$, from which
we immediately derive the general formula through 
$(D^t)^2\psi= t^2\,H^2\psi$. In particular, $H\cdot\psi$ cannot vanish.
\end{proof}
\begin{NB}
Easy examples show that $\psi$ might not be an eigenspinor of $D^t$
itself, since not all constant spinors are eigenspinors of $H$. For
the canonical connection, $\nabla^0T^0=0$ implies that the space of
parallel spinors is invariant under $T^0$, hence there exists a basis
of the space of parallel spinors consisting of eigenspinors.
\end{NB}
\subsection{Vanishing theorems}
This section is devoted to non-existence theorems for solutions
in certain geometric configurations. It allows us to draw quite a
precise picture of what a promising naturally reductive metric 
should look like. First, the Kostant-Parthasarathy formula yields
that we should be interested in precisely those metrics where
$\Omega_{\g}$ is \emph{not} non negative.
\begin{thm}\label{eq-3-4-no-sol}
If the operator $\Omega_{\g}$ is non negative and $\nabla^t$ is not the
Levi-Civita connection, there do not exist any
non trivial solutions to the system of equations
 \bdm
 \nabla^t\psi \ = \ 0, \quad T^t\cdot \psi=0\,.
 \edm
\end{thm}
\begin{proof}
If the spinor $\psi$ is $\nabla^t$-parallel,
then it lies in the kernel of $D^t= D^0 + t H$. Since $\nabla^t$ is 
assumed not to be the Levi-Civita connection, $T^t$ does not vanish and
hence $T^t\cdot \psi=0$ implies $H\cdot \psi=0$. Thus $\psi$ is also 
in the kernel of $D^0$. For the Dirac operator to the parameter $t=1/3$, 
we obtain
\bdm
D^{1/3}\psi \ = \ D^0 \psi + \frac{1}{3} H \cdot \psi \ = \ 0 \,,
\edm
which contradicts Corollary \ref{eigenvalue-estimate}.
\end{proof}
\noindent
For the Levi-Civita connection, it is well known that the
existence of a parallel spinor implies vanishing Ricci curvature.
By repetition of the same argument, one sees that this conclusion does
no longer hold for a metric connection with torsion. Rather,
we get restrictions on the algebraic type of the derivatives of the
torsion. 
\begin{prop}
If the canonical connection $\nabla^0$ is Ricci flat and admits a
parallel spinor, then the  exterior derivative of its torsion $T^0$
satisfies $ (X\haken dT^0)\cdot \psi= 0$ for all vectors $X$ in $\m$.
\end{prop}
\begin{proof}
In \cite[Cor. 3.2]{Friedrich&I01}, Friedrich and Ivanov showed
that a spin manifold with some connection $\nabla$ whose torsion $T$
is totally skew symmetric and a $\nabla$-parallel spinor $\psi$ 
satisfies
\bdm
\left[\frac{1}{2} X\haken dT +\nabla_X T\right]\cdot \psi \ =\
\Ric^{\nabla}(X)\cdot\psi\,.
\edm
Since the canonical connection satisfies $\nabla^0 T^0=0$, the
claim follows.
\end{proof}
\noindent
These conditions are independent of the equation $T^0\cdot\psi=0$.
If $dT^0\neq 0$ and the dimension is sufficiently small, it can happen that 
the intersection of all kernels of $X\haken dT^0$ is already empty, thus 
showing the non-existence of solutions. Models with $dT^0=0$ are of particular
interest and are called \emph{closed} in string theory.\\

\noindent
For further investigations  of the Ricci tensor
 \bdm
 \Ric^t(X,Y)\ =\ \sum_i (t-t^2)\lan [X,Z_i]_{\m},[Y,Z_i]_{\m} \ran +
 Q_{\h}([X,Z_i],[Y,Z_i])\,,
 \edm
it is useful to describe it from a more representation theoretical point
of view. Wang and Ziller derived the general formula we shall present
for $t=1/2$ in \cite{Wang&Z85}. Their proof may easily be generalized to the 
case of arbitrary $t$, hence we omit it here. The main idea is to use a more 
elaborate version of the core computation in the proof of 
Lemma \ref{comp-scalar-in-K-P}. Recall that $C_{\h}$ denotes the
(unlifted) Casimir operator of $\h$, i.\,e.,
\bdm
C_{\h}\ =\ - \sum_{i}\ad X_i \ad Y_i\,.
\edm
It defines a symmetric endomorphism $A: \m\x\m\ra\m$ by
$A(X,Y) := \lan C_{\h}X,Y\ran$.
Similarly, we denote by $\beta(X,Y) = - \tr_{\g} \ad X \ad Y$ the
Killing form of the full Lie algebra $\g$. We make no notational
difference between  $\beta$ itself and its restriction to $\m$.
\begin{thm}\label{gen-Ric-WZ}
The endomorphisms $A$ and $\beta$ satisfy the identities
\bdm
A(X,Y)\ =\  \sum_i Q_{\h}([X,Z_i],[Y,Z_i]),\quad
\beta(X,Y)\ =\ \sum_i\lan [X,Z_i]_{\m},[Y,Z_i]_{\m} \ran +2 A(X,Y)\,.
\edm
Thus, the Ricci tensor is given by
 \bdm
 \Ric^t(X,Y)\ =\ (t-t^2)\beta(X,Y)+ (2t^2-2t+1)A(X,Y)\,.
 \edm
\qed
\end{thm}
\noindent
\begin{NB}
We observe that the coefficient of $\beta$ vanishes for $t=0$ and $t=1$,
and is positive between these parameter values, whereas the
coefficient of $A$ is always positive and attains its minimum for
the Levi-Civita connection ($t=1/2$).
\end{NB}
\noindent
The endomorphism $A$ has block diagonal structure, with every block 
corresponding to an irreducible summand of the isotropy representation.
In particular, the block of the trivial representation
vanishes, since its Casimir eigenvalue is zero. Since $\beta$ is
positive definite for $G$ compact, we can deduce:
\begin{prop}\label{fixvectors-Ric-flat}
Assume that $G$ is compact. If the isotropy representation 
$\Ad: H\ra\SO(\m)$ has fixed vectors, only the connections $t=0$ and $t=1$ 
can be Ricci flat.
\qed
\end{prop}
\noindent
Typically, the eigenvalues of $C_{\h}$ are linear functions of some
deformation parameters, hence, they can vanish for some particular
parameter choices \emph{without} belonging to a trivial $\h$-summand of $\m$.
This makes it difficult to make more precise predictions for
the vanishing of the Ricci tensor.
\begin{prop}\label{H-not-simple}
If the canonical connection has vanishing scalar curvature,
$H$ cannot be simple and the metric cannot be normal homogeneous.
\end{prop}
\begin{proof}
The scalar curvature for the canonical connection is 
\bdm
\sum_{i,j}Q_{\h}([Z_i,Z_j],[Z_i,Z_j])\,.
\edm
By Assumption~\ref{assume-transitivity}, not all vectors
$[Z_i,Z_j]_{\h}$ can be zero. Since $Q_{\h}$ is non degenerate,
we conclude that $Q_{\h}$ can be neither positive nor negative
definite. However, on every simple factor of $\h$, $Q_{\h}$ has to be
a multiple of the Killing form; hence $\h$ cannot be simple.
\end{proof}
\noindent
This fact, as elementary as its proof might be, has far reaching
consequences for the geometry of homogeneous models of string theory.
The existence of a parallel spinor severely restricts the holonomy group
of $\nabla$. In fact, it needs to be a subgroup of the isotropy subgroup
of a spinor inside $\SO(n)$, and these subgroups are well-known. By a 
theorem of Wang 
(\cite[Ch.X, Cor. 4.2]{Kobayashi&N2}), the Lie algebra of the holonomy group 
is spanned by
 \bdm
 \m_0+ [\Lambda_{\m}(\m),\m_0]+ [\Lambda_{\m}(\m),[\Lambda_{\m}(\m),\m_0]]
+\ldots\,,
 \edm
where the subspace $\m_0$ is defined as
 \bdm
 \m_0\ =\ \{[\Lambda_{\m}(X),\Lambda_{\m}(Y)] - \Lambda_{\m}([X,Y]_{\m})
 - \ad([X,Y]_{\h}):\ X,Y\in\m\,  \}\,.
 \edm
For the canonical connection and using our assumption that $[\m,\m]_{\h}$
spans all of $\h$, we conclude that its holonomy Lie algebra is
precisely $\h$. For $t\neq 0$, the holonomy  can only increase, hence
we obtain Table~\ref{table-hol} for the maximally possible subgroups 
$H_{\mathrm{max}}$.
\begin{table}
\bdm
\begin{array}{|c|c|c|c|c|}\hline
\dim G/H^{\vphantom{l^l}} & 5 & 6 & 7 & 8\\[1mm] \hline
H_{\mathrm{max}} & \SU(2)^{\vphantom{l^l}} & \SU(3) & G_2 & \Spin(7)\\[1mm] \hline
\end{array}
\edm
\caption{Maximal holonomy groups for the existence of a parallel spinor}
\label{table-hol}
\end{table}
If we restrict our attention to the canonical connection, 
Proposition~\ref{H-not-simple} implies that $H$ cannot be equal to
$H_{\mathrm{max}}$ itself, but rather has to be a non simple subgroup of it.
This excludes many homogeneous spaces that would naturally come
to one's mind. Of course, they might yield models for other connections
than the canonical one, but such an analysis can only be performed on a
case by case basis.


\section{Examples}\label{examples}
\subsection{The Jensen metric on $V_{4,2}$}
The $5$-dimensional Stiefel manifold $V_{4,2}=\SO(4)/\SO(2)$ carries
a one-parameter family of metrics constructed by G.\ Jensen \cite{Jensen75}
with many remarkable properties. Embed $H=\SO(2)$ into $G=\SO(4)$ as the
lower diagonal $2\x 2$ block. Then the Lie algebra $\so(4)$ splits
into $\so(2)\oplus \m$, where $\m$ is given by
\bdm
\m \ =\left\{ \left[ \ba{c|c} {\ba{cr} 0 & -a \\ a & 0 \ea} & -X^t \\ \hline
X & {\ba{cc} 0 & 0 \\ 0 & 0 \ea} \ea\right] =: (a,X)\, : \
a\in \R,\ X\in \M_{2,2}(\R)\,\right\}\,.
\edm
Denote by $\beta(X,Y):= \tr(X^tY)$ the Killing form of $\so(4)$.
Then the Jensen metric on $\m$ to the parameter $s\in\R$ is  defined
by 
\bdm
\lan (a,X), (b,Y)\ran\ =\ \frac{1}{2}\beta(X,Y) + s \beta(a,b)\ =\
\frac{1}{2}\beta(X,Y)+  2s\cdot ab\,.
\edm 
For $s=2/3$, G.\ Jensen proved that this metric is Einstein, and Th.\
Friedrich showed that it admits two Riemannian Killing spinors 
\cite{Friedrich80} and thus realizes the equality case in his 
estimate for the first eigenvalue of the Dirac operator. A more careful
analysis shows that $V_{4,2}$ carries three different contact structures,
one of which is Sasakian, one quasi-Sasakian but not Sasakian, and the
third one has no special name, although special properties. It will become
clear in the discussion that this metric is only naturally reductive
with respect to $G=\SO(4)$ for $s=1/2$. In the following sections,
we shall describe the Jensen metrics on $V_{4,2}$ first from the point of 
view of contact geometry and then from the point of view of naturally
reductive spaces.
%
\subsection{The contact geometry approach}
%
Denote by $E_{ij}$ the standard basis of $\so(4)$. Then the elements
 \bdm
 Z_1\, :=\, E_{13},\ Z_2 \,:=\, E_{14},\ Z_3\,=\,E_{23},\ Z_4\,=\,
 E_{24},\ Z_5\,=\,\frac{1}{\sqrt{2s}}\,E_{12}
 \edm
form an orthonormal base of $\m$. To start with, we compute all 
nonvanishing commutators in $\m$. These are
 \begin{equation}\tag{$*$}
 \ba{ccccccccc}
 [Z_1,Z_3]_{\m}\!& = &\!\sqrt{2s}\,Z_5, & 
 [Z_1,Z_5]_{\m}\!& = &\! - \displaystyle{\frac{1}{\sqrt{2s}}}
 \,Z_3, & [Z_2,Z_4]_{\m}\!& = &\!\sqrt{2s}\,Z_5, \\[4mm]
 [Z_2,Z_5]_{\m}\! & = &\!-\displaystyle{\frac{1}{\sqrt{2s}}}\,Z_4, & 
 [Z_3,Z_5]_{\m}\!& = &\! \displaystyle{\frac{1}{\sqrt{2s}}}\,Z_1, & 
 [Z_4,Z_5]_{\m}\!& = &\! \displaystyle{\frac{1}{\sqrt{2s}}}\,Z_2.
 \ea\end{equation}
Notice that all these commutators have no $\h$-contribution.
Identifying $\m$ with $\R^5$ via the chosen basis, the isotropy 
representation of an element $g(\theta)=\left[\ba{cc}\cos\theta& -\sin\theta\\
\sin\theta&\cos\theta \ea\right]\in H=\SO(2)$ may be written as follows:
 \bdm
 \Ad g(\theta)\ =\  \left[\ba{ccccc}
 \cos\theta & -\sin\theta & 0 & 0 & 0 \\
 \sin\theta & \cos\theta & 0 & 0 & 0 \\
 0 & 0 & \cos\theta & -\sin\theta & 0 \\
 0 & 0 & \sin\theta & \cos\theta & 0 \\
 0 & 0 & 0 & 0 & 1   
 \ea\right]\,.
 \edm
In particular, $Z_5$ is invariant under the isotropy action.
As in \cite{Friedrich80}, we use a suitable  basis $\psi_1,\ldots,\psi_4$
for the $4$-dimensional spinor 
representation $\kappa: \Spin(\R^5)\ra \GL(\Delta_5)$. One  derives 
the expression for the lift of the isotropy representation,
 \bdm
 \kappa\big(\Adtilde g(\theta)\big)\ =\ \left[\ba{cccc}
 \cos\theta+i\sin\theta & 0 & 0 & 0 \\
 0 & \cos\theta-i\sin\theta & 0 & 0 \\
 0 & 0& 1 & 0\\ 0 & 0 & 0 & 1 \ea\right]\,.
 \edm
Thus, the elements $\psi_3$ and $\psi_4$ define sections of 
the spinor bundle $S=G\x_{\kappa(\Adtilde)}\Delta_5$ if viewed as 
constant maps $G\ra\Delta_5$. In fact, for $s=2/3$, 
$\psi^{\pm}:= \psi_3 \mp i\psi_3$
are exactly the Riemannian Killing spinors from \cite{Friedrich80} as
we will see below. The sections induced by  $\psi_1$ and $\psi_2$
are not constant and thus more difficult to treat. We will not consider
them in our discussion.
In \cite[Prop. 3]{Jensen75}, the author computed the
map $\Lambda_{\m}^{\LC}: \m\cong\R^5\ra\so(5)$ (see Wang's Theorem in
Section~\ref{fam-conn}) defining the Levi-Civita connection:
 \bdm
 \Lambda_{\m}^{\LC}(Z_{\alpha})Z_{\beta} = 
 \frac{1}{2}[Z_{\alpha},Z_{\beta}],\
 \Lambda_{\m}^{\LC}(Z_5)Z_{\alpha} = (1-s) [Z_5,Z_{\alpha}],\
 \Lambda_{\m}^{\LC}(Z_{\alpha})Z_5 = s [Z_{\alpha},Z_5]\quad \text{for }
 \alpha, \beta= 1,\ldots,4.
  \edm
Indeed, one easily verifies that this is the unique map $\Lambda_{\m}$
verifying the conditions
 \bdm
 \lan \Lambda_{\m}(X)Y,Z\ran + \lan Y,\Lambda_{\m}(X)Z\ran\ =\ 0\ \text{ and }\
 \Lambda_{\m}(X)Y - \Lambda_{\m}(Y)X \ =\ [X,Y]_{\m}\,.
 \edm
Thus, one sees that for $s\neq 1/2$, $\Lambda_{\m}(X)Y$ is not globally
proportional to the commutator $[X,Y]_{\m}$, and both
$\lan \Lambda_{\m}(X)Y,Z\ran$ and  $- \lan [X,Y]_{\m},Z\ran$ (the torsion
of the canonical connection) fail to define a $3$-form: The first is not
skew symmetric in $X$ and $Y$, the second is not skew symmetric in $X$ and
$Z$. In any case,
by using the commutator relations $(*)$, the Levi-Civita connection can
be identified with an endomorphism of $\R^5$ as follows:
 \bdm
 \Lambda_{\m}^{\LC}(Z_1)\,=\, \sqrt{\frac{s}{2}}E_{35},\quad 
 \Lambda_{\m}^{\LC}(Z_2)\,=\, \sqrt{\frac{s}{2}}E_{45},\quad 
 \Lambda_{\m}^{\LC}(Z_3)\,=\, - \sqrt{\frac{s}{2}}E_{15},\quad 
 \Lambda_{\m}^{\LC}(Z_4)\,=\, - \sqrt{\frac{s}{2}}E_{25},
\edm\bdm
 \Lambda_{\m}^{\LC}(Z_5)\,=\, \frac{1-s}{\sqrt{2s}}(E_{13}+E_{24})\,.
 \edm
The lift into the spin representation  yields a global factor $1/2$
and replaces $E_{ij}$ by $Z_i\hut Z_j$. By setting
 \bdm
 \tilde{T}\ :=\ (Z_1\hut Z_3+Z_2\hut Z_4)\hut Z_5 ,
 \edm
the Levi-Civita connection may be rewritten in a unified way as
 \begin{equation}\label{LC-uniform}
 \tilde{\Lambda}_{\m}^{\LC}(Z_5)\,=\, \frac{1}{4}\frac{2(1-s)}{\sqrt{2s}}
 (Z_5\haken \tilde{T}),\quad 
 \tilde{\Lambda}_{\m}^{\LC}(Z_{\alpha})\,=\, \frac{1}{4}\sqrt{2s}
 (Z_{\alpha}\haken \tilde{T}) \ \text{ for }\alpha=1,\ldots,4\,.
 \end{equation}
Now we discuss the three different metric almost contact structures existing 
on $V_{4,2}$. The space $\m$ has a preferred direction, namely
$\xi= Z_5$, which is fixed under the isotropy representation. 
Denote its dual $1$-form, $\eta(X)=\lan Z_5, X\ran$ by $\eta$. 
The following operators
 \bdm
 \varphi_{S}\,=\, \left[\ba{ccccc}0 & 0 & 1 &0&0\\ 0&0&0&1&0\\
 -1 & 0 & 0&0&0\\ 0& -1& 0&0&0\\ 0&0&0&0&0 \ea \right]\,,\
 \varphi_{qS}\,=\, \left[\ba{ccccc}0 & 1 & 0 &0&0\\ -1&0&0&0&0\\
 0 & 0 & 0&1&0\\ 0& 0& -1&0&0\\ 0&0&0&0&0 \ea \right]\,,\
 \varphi_{*} \,=\, \left[\ba{ccccc}0 & 1 & 0 &0&0\\ -1&0&0&0&0\\
 0 & 0 & 0&-1&0\\ 0& 0& 1&0&0\\ 0&0&0&0&0 \ea \right]
 \edm
intertwine the isotropy representation, and thus define compatible
complex structures on the linear span of $Z_1,\ldots,Z_4$. Then one checks
for all three choices for $\varphi$  that the compatibility conditions 
defining a metric almost contact structure hold:
 \bdm
 \varphi^2\,=\, -\mathrm{Id} + \eta\ox\xi,\quad 
 \lan\varphi(X),\varphi(Y)\ran\,=\, \lan X,Y\ran - \eta(X)\cdot\eta(Y),\quad 
 \varphi(\xi)\,=\, 0\,.
 \edm
The fundamental
form of the structure is defined by $F(X,Y)= \lan X,\varphi(Y) \ran$, thus
yielding
 \bdm
 F_{S}\ =\ Z_1\hut Z_3 + Z_2\hut Z_4,\quad 
 F_{qS}\ =\ Z_1\hut Z_2 + Z_3\hut Z_4,\quad 
 F_{*}\ =\ Z_1\hut Z_2 - Z_3\hut Z_4\,,
 \edm
respectively. Since $Z_5$ is constant under the isotropy action, its exterior 
derivative may be computed using the general formula as stated at the
beginning of the proof of Lemma \ref{ext-diff},
 \bdm
 d \omega^1(X_0,X_1)\ =\ X_0(\omega^1(X_1)) -  X_1(\omega^1(X_0))
 - \omega^1([X_0,X_1])\,.
 \edm
For the constant vector field $Z_5$, we thus obtain
$d Z_5(Z_i,Z_j) = - \lan Z_5, [Z_i,Z_j]\ran$. Applying again the commutator 
relations implies 
 \bdm
 dZ_5\  =\ - \sqrt{2s}\,(Z_1\hut Z_3+Z_2\hut Z_4)\,.
 \edm
In particular, $dZ_5$ is
proportional to $F_{S}$, turning it into a Sasaki structure (up to rescaling)
and  implying immediately $d F_{S}=0$.  For the other two structures,
remark that $Z_1\hut Z_2$ and $Z_3\hut Z_4$ are also invariant forms
under the isotropy action, thus their exterior differential may be
computed in a similar way. One gets that $dF_{qS}=0$, turning it
into a non Sasakian quasi-Sasakian structure, and 
$dF_*$ is proportional to $Z_2\hut Z_3\hut Z_5$, which implies 
$dF^{\varphi^*}_*=0$. We can then compute the Nijenhuis tensor
 \bdm
 N(X,Y)\ :=\ [\vphi(X), \vphi(Y)] + \vphi^2([X,Y]) - \vphi([\vphi(X),Y])
 - \vphi([X,\vphi(Y)]) + d\eta(X,Y)\cdot \xi
 \edm
and see that it vanishes for all three metric almost contact structures. 
By \cite[Thm. 8.2]{Friedrich&I01}, the Stiefel manifold $V_{4,2}$ admits 
a unique almost contact connection $\nabla$ with torsion
 \bdm
 T\,=\, \eta\hut d\eta \, =\,- \sqrt{2s}\,(Z_1\hut Z_3+Z_2\hut Z_4)\hut Z_5\,.
 \edm
Next we discuss the existence of spinors that
are parallel with respect to the connection $\nabla$
as well as the existence of Killing spinors, since we consider the analogy and
differences to the previous case to be instructive.
\begin{thm}
\begin{enumerate}
\item[]
\item The constant spinors are parallel with respect to
the contact connection $\nabla$  if and only if $s=1/2$;
\item The constant spinors $\psi^{\pm}$ are Riemannian Killing spinors 
if and only if $s=2/3$.
\end{enumerate}
\end{thm}
\begin{proof}
In equation~(\ref{LC-uniform}), we gave the general formula
for the Levi-Civita connection in direction $Z_i$ as the inner product of
$Z_i$ and the $3$-form $\tilde{T}$. If a constant spinor $\psi$ is to
be parallel with respect to  $\nabla$,
 \bdm
 0\ =\ \nabla_X \psi\ =\ (\tilde{\Lambda}^{\LC}_{\m}(X) + 
 \frac{1}{4}X\haken T ) \psi\,,
 \edm
then the coefficients of $\tilde{\Lambda}^{\LC}_{\m}$ as in 
equation~(\ref{LC-uniform}) have to be equal for all $Z_i$, hence,
$2(1-s)/\sqrt{2s}=\sqrt{2s}$, which means that $s=1/2$. For this value,
the combination $\tilde{\Lambda}^{\LC}(X) + \frac{1}{4}X\haken T$
vanishes, so both constant spinors are parallel indeed. For the
discussion of Riemannian Killing spinors, we use the following realization of 
the spin representation:
\bdm\ba{ccc}
e_1=\left[\ba{cccc} 0&0&0&i\\0&0&i&0\\0&i&0&0\\i&0&0&0\ea\right], & 
e_2=\left[\ba{cccc} 0&0&0&-1\\0&0&1&0\\0&-1&0&0\\1&0&0&0\ea\right], &
e_3=\left[\ba{cccc} 0&0&-i&0\\0&0&0&i\\-i&0&0&0\\0&i&0&0\ea\right], \\
& & \\
e_4=\left[\ba{cccc} 0&0&1&0\\0&0&0&1\\-1&0&0&0\\0&-1&0&0\ea\right], & 
e_5=\left[\ba{cccc} i&0&0&0\\0&i&0&0\\0&0&-i&0\\0&0&0&-i\ea\right]. &
\ea\edm
Then one checks that 
\bdm
(Z_5\haken \tilde{T})\cdot \psi^{\pm}\ =\ \pm 2 Z_5\cdot \psi^{\pm},\quad
(Z_{\alpha}\haken \tilde{T})\cdot \psi^{\pm}\ =\ \pm  Z_{\alpha}\cdot 
\psi^{\pm} \text{ for }\alpha=1,\ldots,4.
\edm
Looking at $Z_5$, we conclude that the Killing equation 
$\nabla^{\LC}_X\psi=\mu\, X\cdot\psi$ implies that the coefficients  
in equation~(\ref{LC-uniform}) have to satisfy
$2(1-s)/\sqrt{2s}=\sqrt{2s}/2$. The solution is now $s=2/3$, and one
checks that $\psi^{\pm}$ are Killing spinors indeed.
\end{proof}
%
\subsection{The naturally reductive space approach}
%
We would like to interpret the metric $\lan\,,\,\ran$ as a naturally reductive
metric with respect to some other group $\bar{G}$, and the connection with 
the torsion
 \bdm
 T \, =\,- \sqrt{2s}\,(Z_1\hut Z_3+Z_2\hut Z_4)\hut Z_5
 \edm
as its canonical  connection. So write $M=\bar{G}/\bar{H}$ with the
Lie algebra decomposition $\bar{\g}=\bar{\h}\oplus\bar{\m}$, and
assume that the original isotropy representation is a subrepresentation of 
the new isotropy representation, i.\,e., the action of $\h\subset\bar{\h}$ 
on $\m\cong\bar{\m}$ remains unchanged. This point of view necessarily 
enlarges the holonomy group $H$ already for dimensional reasons.
In fact, we can deduce a lot of information about the new isotropy
representation from the formula for $T$. In Remark~\ref{torsion-lie-alg}, 
we explained the relation between $\m$-commutators
and the torsion. For example, the formula above implies
\bdm
[Z_1,Z_3]_{\bar{\m}}\ =\ \sqrt{2s}Z_5,\quad
[Z_4,Z_5]_{\bar{\m}}\ =\ \sqrt{2s}Z_2,\quad
[Z_1,Z_4]_{\bar{\m}}\ =\ [Z_3,Z_4]_{\bar{\m}}\ =\ 0 \,.
\edm
Then we can compute
\bdm
\Jac_{\bar{\m}}(Z_1,Z_3,Z_4) \ =\ 2s\,Z_2\,.
\edm
On the other hand, 
\bdm
\Jac_{\bar{\h}}(Z_1,Z_3,Z_4)\ =\ - Z_2+ [Z_4,[Z_1,Z_3]_{\bar{\h}}]+
[Z_3,[Z_4,Z_1]_{\bar{\h}}]\ \stackrel{!}{=} - \Jac_{\bar{\m}}(Z_1,Z_3,Z_4)\,.
\edm
Thus, there must be two elements $H_1:=[Z_1,Z_3]_{\bar{\h}}$ and 
$H_2:=[Z_4,Z_1]_{\bar{\h}}$ in $\bar{\h}$, not both zero, such that
 \bdm
 [H_1,Z_4]+[H_2,Z_3]\ =\ (2s-1) Z_2\,.
 \edm
By some more careful analysis, one obtains  $H_2=0$,
$H_1=[Z_2,Z_4]_{\bar{\h}}$ and the action of $H_1$ on the other
vectors $Z_i$. The systematic description of $\lan\,,\,\ran$ as a
naturally reductive metric can be given using a deformation
construction due to Chavel and Ziller (\cite{Chavel70}, \cite{Ziller77}). 
It is based on the remark that for $s=1/2$, $\m$ splits into an 
orthogonal direct sum 
of $\m_1:=\{(0,X)\}$ and $\m_2:=\{(a,0)\}$ such that
 \bdm
 [\h,\m_2]\ =\ 0 \text{ and }[\m_2,\m_2]\subset \m_2\,.
 \edm
Let $M_2\subset G$ be the subgroup of $G$ with Lie algebra $\m_2$, and
set $\bar{G}=G\x M_2$, $\bar{H}=H\x M_2$. An element $(k,m)$ of $\bar{G}$
acts on $M=G/H$ by $(k,m)gH=kgHm^{-1}$, and then $\bar{H}$ can  indeed
be identified with the isotropy group of this action. We endow
$\bar{\g}=\g\oplus \m_2$ with the direct sum Lie algebra structure. 
The trick is now to choose a realization of $\bar{\m}$ that depends
on the deformation parameter $s$ of the metric. Writing all elements
of $\bar{\g}$ as $4$-tuples $(H,U,X,Y)$ with $H\in\h$, $U\in\m_1$
and $X,Y\in\m_2$, we can realize the Lie algebra of $\bar{H}$ as
 \bdm
 \bar{\h}\ = \ \{(H,0,X,X)\subset \bar{\g}:\ H\in \h,\, X\in\m_2\}
 \edm
and choose 
 \bdm
 \m\ =\ \{(0,X, 2s\,Y,(2s-1)Y) :\ X\in\m_1, Y\in\m_2\}
 \edm
as an orthogonal complement.
Here, $(0,0,2s\,Y,(2s-1)Y)$ will be identified with $Y\in\m_2$. 
Since $\m_2$ is abelian in this example, the Lie algebra structure
of $\bar{\g}$ is particularly simple. $\bar{\h}$ is a Lie
algebra with commutator
 \bdm
 [(H,0,X,X),(H',0,X',X')]\ =\ ([H,H'],0,0,0)\,,
 \edm
the full isotropy representation is
 \bdm
 [(H,0,X,X),(0,U,2sY,(2s-1)Y)]\ =\ (0,[H+X,U],0,0)
 \edm
and the commutator of two elements in $\bar{\m}$ splits into
its $\bar{\h}$ and $\bar{\m}$ part as follows:
 \bea[*]
 [(0,U,2s\,X,(2s-1)X),(0,V,2s\,Y,(2s-1)Y)]& =&
 ([U,V]_{\h},0,-(2s-1)[U,V]_{\m_2},-(2s-1)[U,V]_{\m_2} )\\
 & &\hspace{-2.3cm} +
 (0,[U,V]_{\m_1}+2s([U,Y]+[X,V]), 2s[U,V]_{\m_2},(2s-1)[U,V]_{\m_2} )\,.
 \eea[*]
With these choices for $\bar{\h}$ and $\bar{\m}$, the metric
$\lan\,,\,\ran$ is naturally reductive with respect to $\bar{G}$, 
the torsion of its canonical connection is precisely $T$ and
the Ricci tensor is given by
 \bdm
 \Ric^0\ =\ 2(1-s)\diag(1,1,1,1,0)\,.
 \edm
For $s=1$, the canonical connection is thus Ricci flat, and by
Proposition~\ref{fixvectors-Ric-flat}, we know that no other
connection can have this property. However, the holonomy 
$\bar{H}\cong \SO(2)\x\SO(2)$ is too large to admit parallel spinors.
For $s=1/2$, we have two parallel spinors for the canonical connection
as seen in the preceding section, but the Ricci curvature does not vanish. 
In this 
case, one can ask the question whether some other connection of
the family $\nabla^t$ admits parallel spinors. But using 
Wang's Theorem (\cite[Ch.X, Cor. 4.2]{Kobayashi&N2}) for computing the
holonomy, one sees that  $\nabla^t$ has full holonomy $\SO(\m)$ for
$t\neq 0$, excluding again the existence of parallel spinors.

We close this section with a look at the eigenvalue estimate
for $(D^{1/3})^2$. Since the extension of $H$ is by the abelian group
$\SO(2)$, the Casimir operator $\Omega_{\g}$ is non negative by
Lemma~\ref{non-negative} and Corollary~\ref{eigenvalue-estimate} can be
applied. We compute the scalar in the general Kostant-Parthasarathy
formula (Theorem~\ref{K-P-1})
\bdm
\frac{1}{8}\sum_{i,j} Q_{\h}([Z_i,Z_j],[Z_i,Z_j])
 +\frac{3}{8}t^2 \sum_{i,j} Q_{\m}([Z_i,Z_j],[Z_i,Z_j])\, =\,
\frac{1}{8}\cdot 8(1-s)+ \frac{3}{8}t^2\cdot 24 s\, =\,
1+ (9t^2-1)s
\edm
and see that it is independent of the deformation parameter $s$
precisely for the Kostant connection $t=1/3$. If $s\neq 1/2$, there exist
no constant spinors and hence Corollary~\ref{eigenvalue-estimate} is a strict 
inequality,
 \bdm
 (\lambda^{1/3})^2\ > \ 1\,.
 \edm
For $s=1/2$, there
exists a constant spinor $\psi$ and it satisfies by 
Theorem~\ref{constant-fields}
\bdm
(D^t)^2\psi \ =\ 9t^2\cdot 1\cdot \psi \ =\ 9t^2\psi \,.
\edm
Unfortunately, we have been unable to relate  this bound with the
infimum of the spectrum of $(D^t)^2$ for other values of $t$.
In particular, it seems to be difficult to deduce from 
Corollary~\ref{eigenvalue-estimate} any information about the Riemannian
Dirac spectrum.

\providecommand{\bysame}{\leavevmode\hbox to3em{\hrulefill}\thinspace}

\end{document}